\documentclass[a4paper,12pt]{article}
\pdfoutput=1
\usepackage{etex}

\usepackage[inline]{asymptote}

\usepackage{amssymb}
\usepackage{amsmath}
\usepackage{amsthm}
\usepackage{bbm}

\usepackage{lineno}
\usepackage[cm]{fullpage}
\usepackage{graphicx}
\usepackage{wrapfig}
\usepackage{subfig}
\usepackage{float}
\usepackage{booktabs}
\usepackage{caption}

\allowdisplaybreaks

\makeatletter
\newsavebox{\sfe@box}
\newenvironment{subfloatenv}[1]{%
\def\sfe@caption{#1}%
\setbox\sfe@box\hbox\bgroup\color@setgroup}%
{\color@endgroup\egroup\subfloat[\sfe@caption]%
{\usebox{\sfe@box}}}
\makeatother

\theoremstyle{definition}

\newcommand{\norm}[1]{\ensuremath{\lVert#1\rVert}}

\newcommand{\reverse}[1]{\ensuremath{\widetilde{#1}}}

\newcommand{\R}[1]{\ensuremath{\mathbb{R}^{#1}}}

\newcommand{\El}[1]{\ensuremath{{E}_{#1}}}
\newcommand{\Hy}[1]{\ensuremath{{H}_{#1}}}

\newcommand{\tb}[1]{\ensuremath{\textbf{#1}}}
\newcommand{\mb}[1]{\ensuremath{\boldsymbol{#1}}}

\newcommand{\e}{\tb{e}} % basis vectors
\newcommand{\I}{\tb{I}} % pseudoscalar
\newcommand{\J}{\ensuremath{J}} % duality transform; isomorphism from R^n to R^n*
\title{Clifford algebra and the projective model of Hyperbolic spaces}
\author{Andrey Sokolov}

\begin{document}
\maketitle

\begin{abstract}
I apply the algebraic framework developed in \cite{gunn2011geometry} to study geometry of hyperbolic spaces in 1, 2, and 3 dimensions.
The background material on projectivised Clifford algebras and their application to Cayley-Klein geometries
is described in \cite{sokolov2013clifford}. 
\end{abstract}

\tableofcontents

\section{Hyperbolic line \Hy{1}}
For a point \(\tb{a}=d\e_0+a\e_1\), I compute \(\tb{a}\reverse{\tb{a}}=\tb{a}^2=a^2-d^2\)
and \(\norm{\tb{a}}=\sqrt{|a^2-d^2|}\).
A point \(\tb{a}\) is proper if \(\tb{a}^2>0\), which is equivalent to \(|a|>|d|\).
Therefore, \(\tb{a}=-x\e_0+\e_1\) is proper only if \(|x|<1\).
Positively oriented normalised proper points can be parameterised by \(\phi\) as follows:
\begin{equation}\label{parameterise Hy1}
\tb{a}=-\e_0\sinh{\phi}+\e_1\cosh\phi.
\end{equation}
As a vector, \(\tb{a}\) dually represents a point in \Hy{1} at \(x=\tanh\phi\).
This parametrisation does not cover negatively oriented points.
Since \(\tb{a}=\e_1(\cosh\phi+\e_{01}\sinh\phi)=\e_1e^{\phi\e_{01}}\), I can write
\begin{equation}
\tb{a}=\e_1e^{\phi\e_{01}}
\end{equation}
for any normalised proper point with the positive orientation.

In \Hy{1}, the distance \(r\) between normalised proper points \(\tb{a}\) and \(\tb{b}\) is defined by 
\begin{equation}
\sinh{r}=|\tb{a}\vee\tb{b}|.
\end{equation}
The distance also satisfies \(\cosh{r}=|\tb{a}\cdot\tb{b}|\), 
since \(|\tb{a}\cdot\tb{b}|^2-|\tb{a}\vee\tb{b}|^2=1\) for any normalised proper point.
For \(\tb{a}=-\e_0\sinh{\phi}+\e_1\cosh\phi\) and \(\tb{b}=-\e_0\sinh{\theta}+\e_1\cosh\theta\),
I get \(\tb{a}\cdot\tb{b}=\cosh\phi\cosh\theta-\sinh\phi\sinh\theta=\cosh(\phi-\theta)\) and, therefore,
\begin{equation}
r=|\phi-\theta|.
\end{equation}
It follows that the point \(\tb{n}_+=-\e_0+\e_1\) is at the positive infinity, \(r=+\infty\),
and the point \(\tb{n}_-=\e_0+\e_1\) is at the negative infinity,  \(r=-\infty\)
(both \(\tb{n}_+\) and \(\tb{n}_-\) are null).
The distance to the other points is undefined, e.g.\ the distance from the origin to \(\e_0\) is undefined.
The metric component of the space \Hy{1} resides within the interval \((-1,1)\).
Both \(x=+1\) and \(x=-1\) are at the infinite distance from the origin.
The points outside of \((-1,1)\) exists in \Hy{1} in a non-metric sense only, since they are not reachable by translation from 
anywhere in the interval \((-1,1)\). 
So, in the metric sense, the space \Hy{1} is open and non-periodic, with two distinct infinities at \(x=+1\) and \(x=-1\).

\begin{figure}[h]
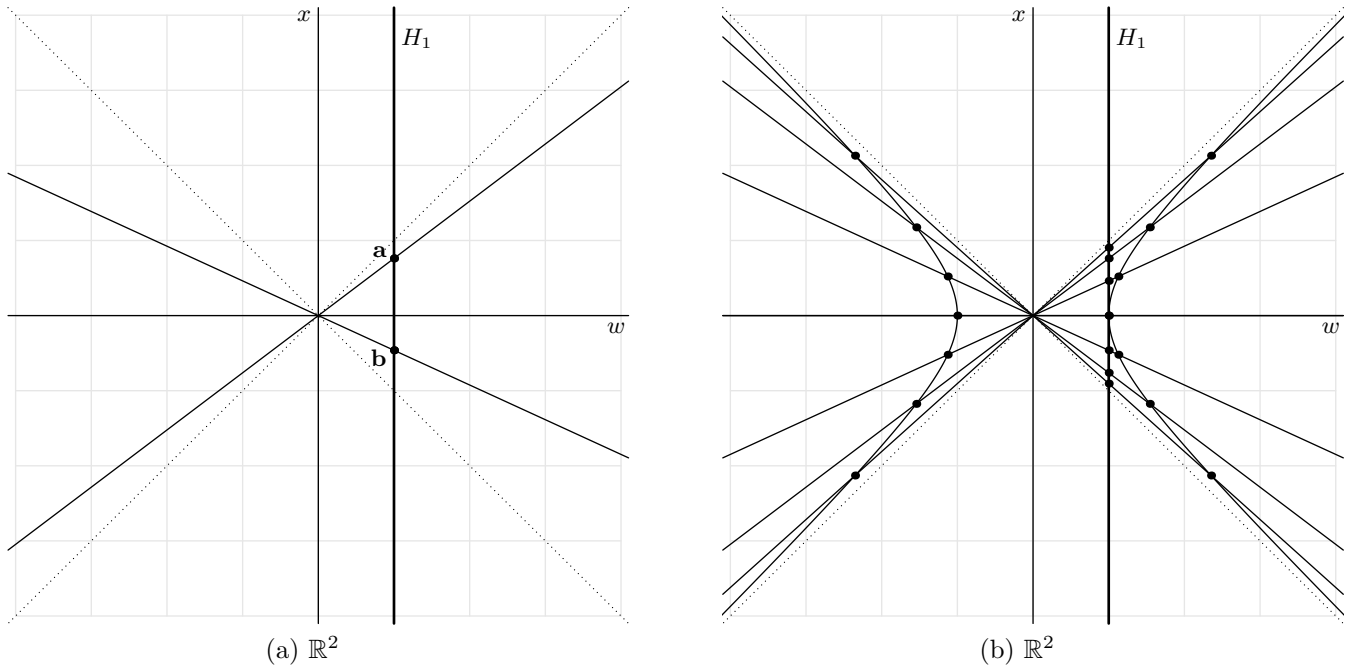

\begin{subfloatenv}{\R{2}}
\begin{asy}
import Figure2D;
Figure f = Figure(xaxis_name="$w$",yaxis_name="$x$");

f.null_lines();
f.line(e_1-e_0,draw_orientation=false,label="$H_1$",position=0.05,align=(1,0),pen=currentpen+1);

real alpha = 1;
f.line(Line(0,-sinh(alpha),cosh(alpha)),draw_orientation=false);
dot(Label("$\textbf{a}$",align=(-1,0.5)),(1,tanh(alpha)));
//f.arc1(O,Point(0,0,1),Point(0,3,1),direction=CW,radius=0.4,label="$\theta$",align=(0.5,0));
dot((1,tanh(alpha)));

real beta = -0.5;
f.line(Line(0,-sinh(beta),cosh(beta)),draw_orientation=false);
dot(Label("$\textbf{b}$",align=(-1,-0.5)),(1,tanh(beta)));
//f.arc1(O,Point(0,0,1),Point(0,-2,1),radius=0.35,label="$\alpha$",align=(0.5,0));
dot((1,tanh(beta)));

//f.arc(O,Point(0,-2,1),Point(0,-3,-1),radius=1.6,label="$r_{ac}$");

\end{asy}
\end{subfloatenv}\hfill%
\begin{subfloatenv}{\R{2}}
\begin{asy}
import Figure2D;
Figure f = Figure(xaxis_name="$w$",yaxis_name="$x$");

real alpha_max = 2.5;
real x(real alpha) { return cosh(alpha); }
real y(real alpha) { return sinh(alpha); }
path p = graph(x, y, -alpha_max, alpha_max, operator ..);
draw(p);

real x(real alpha) { return -cosh(alpha); }
real y(real alpha) { return sinh(alpha); }
path p = graph(x, y, -alpha_max, alpha_max, operator ..);
draw(p);

f.crop();
f.null_lines();

f.line(e_1-e_0,draw_orientation=false,label="$H_1$",position=0.05,align=(1,0),pen=currentpen+1);

int n=3;
for(int i=-n; i<=n; ++i) { real alpha = i/2.0; write(alpha); MV a = Line(0,-sinh(alpha),cosh(alpha)); dot((1,tanh(alpha))); dot((cosh(alpha),sinh(alpha))); dot((-cosh(alpha),sinh(alpha))); f.line(a,draw_orientation=false);}

\end{asy}
\end{subfloatenv}
\caption{The distance between points in \Hy{1}.}
\label{distance in Hy1}
\end{figure}

The points \(\tb{a}\) and \(\tb{b}\) shown in Figure~\ref{distance in Hy1}(a) are
parameterised with \(\phi=1\) and \(\theta=-\tfrac{1}{2}\), respectively.
So, the distance between them evaluates to \(r=|1-(-\tfrac{1}{2})|=\tfrac{3}{2}\).
The subspaces dually represented by the null vectors \(\tb{n}_+\) and \(\tb{n}_-\) are shown with dotted lines.

It is sometimes useful to visualise \Hy{1} as a pair of hyperbolic curves defined by
\begin{equation}
w^2-x^2=1
\end{equation}
and shown in Figure~\ref{distance in Hy1}(b), instead of a straight line at \(w=1\).
Antipodal points on the curves are identified, 
but it is sufficient to consider just one branch corresponding to, say, \(w\ge1\).
In contrast with the circular model of \El{1},
the relationship between the distance in \Hy{1} and the arc length on the hyperbolic curves is not uniform.
For example, a series of points equidistant in \Hy{1} is shown in Figure~\ref{distance in Hy1}(b);
the points are parameterised with \(\phi=-1.5,-1,-0.5,0,0.5,1,1.5\) (the distance between the adjacent points is 0.5).
Clearly, these points are not separated by segments of equal length.
In fact, it is the areas enclosed by the hyperbolic curves and the relevant straight lines
that are equal.
Moreover, the space \Hy{1} is characterised by a constant negative curvature of \(-1\), whereas
a hyperbola has a variable positive curvature.
The utility of the hyperbolic curves is that they provide a model of \Hy{1} in which it is clearly unbound from Euclidean point of view.

Hyperbolic space is superficially similar to Minkowski space, 
since both spaces make use of the hyperbolic measure (related to the use of hyperbolic functions).
However, in Minokwski space it is the angular measure that is hyperbolic
and in hyperbolic space it is the distance measure that is hyperbolic.
Moreover, Minkowski space is flat (zero curvature) and the distance measure is degenerate there,
whereas hyperbolic space has a constant negative curvature, which gives it the hyperbolic distance measure.

The standard definition of spinors applies in \Hy{1}.  Any spinor can be written as either \(S=e^{\phi\e_{01}}\)
or  \(-e^{\phi\e_{01}}\) for some \(\phi\in\R{}\).
Proper motions consist of translations and are generated by the action \(S\tb{a}S^{-1}\).
Note that the action of \(-S\) is exactly the same as
that of \(S\), so it is sufficient to consider spinors \(S=e^{\phi\e_{01}}\).

As usual, the translation of a point \(\tb{a}\) by \(\lambda\) is given by \(T\tb{a}T^{-1}\), 
where \(T=e^{-\tfrac{1}{2}\lambda\e_{01}}\),
and if \(\tb{a}\) is parameterised with \(\phi\), then
\(\tb{a}'=T\tb{a}T^{-1}\) can be parameterised with \(\phi'=\phi+\lambda\).
So, \(\e_1e^{\phi\e_{01}}\) can be interpreted as the translation of the origin \(\e_1\) by \(\phi\).

The top-down reflection of \(\tb{a}\) in \(\tb{b}\) is given by \(-\tb{b}\tb{a}\tb{b}^{-1}\),
and if \(\tb{a}\) and \(\tb{b}\) are parameterised with \(\phi\) and \(\theta\), respectively,
then \(\tb{a}'=-\tb{b}\tb{a}\tb{b}^{-1}=-\e_1e^{(\theta-\phi+\theta)\e_{01}}\).
If \(\tb{a}\) is positively oriented, then \(\tb{a}'\) is negatively oriented, and vice versa.

If \(\tb{a}\) is a proper point at \(x\in(-1,1)\), then its polar point \(\tb{a}\e_{01}\) is at \(1/x\) and, therefore, is not proper.
The points \(\tb{n}_+\) and \(\tb{n}_-\) are polar to themselves (up to a change in orientation in the case of \(\tb{n}_-\)),
i.e.\ \(\tb{n}_+\e_{01}=\tb{n}_+\) and \(\tb{n}_-\e_{01}=-\tb{n}_-\).

The projection of \(\tb{a}\) on \(\tb{b}\) is given by \((\tb{a}\cdot\tb{b})\tb{b}^{-1}\),
which  yields \(\tb{b}\cosh(\phi-\theta)\) if \(\tb{a}\) and \(\tb{b}\) are proper and normalised
(\(\tb{a}\) on \(\tb{b}\) are parameterised by \(\phi\) and \(\theta\), respectively).
The rejection of  \(\tb{a}\) by \(\tb{b}\) is given by \((\tb{a}\wedge\tb{b})\tb{b}^{-1}\),
which  yields \(\tb{b}\e_{01}\sinh(\phi-\theta)\) for the normalised proper points.
As in elliptic space, the rejection by \(\tb{b}\) is equivalent to the projection on the polar point of \(\tb{b}\).
It is not possible to project on or reject by \(\tb{n}_+\) and \(\tb{n}_-\) since they are null and not invertible.

\section{Hyperbolic plane \Hy{2}}
For a line \(\tb{a}=d\e_0+a\e_1+b\e_2\), I have \(\tb{a}\reverse{\tb{a}}=\tb{a}^2=a^2+b^2-d^2\) and
\(\norm{\tb{a}}=\sqrt{|a^2+b^2-d^2|}\). 
The line is proper if \(\tb{a}^2>0\), null if \(\tb{a}^2=0\), and improper if \(\tb{a}^2<0\).
For a point \(\tb{P}=w\e_{12}+x\e_{20}+y\e_{01}\), I have \(\tb{P}\reverse{\tb{P}}=-\tb{P}^2=w^2-x^2-y^2\) and
\(\norm{\tb{P}}=\sqrt{|w^2-x^2-y^2|}\).
The point is proper if \(\tb{P}^2<0\), null if \(\tb{P}^2=0\), and
improper if \(\tb{P}^2>0\).
Recall that \(\tb{P}=\e_{12}+x\e_{20}+y\e_{01}\) represents a counterclockwise point in \Hy{2} at \((x,y)\).
It is proper if \((x,y)\) is located within the unit circle\footnote{Note that this descriptive terminology is used in a non-metric sense here.} 
defined by
\begin{equation}
x^2+y^2=1,
\end{equation}
improper if it is outside this circle, and null if it lies on the circle.
The same applies to clockwise points and points with arbitrary weight.
A line is proper if it intersects the unit circle in two points, or alternatively if the line passes through a proper point.
It is null if it is tangent to the unit circle, and improper otherwise.

%Any normalised proper counterclockwise point in \Hy{2} can be written as
%\begin{equation}
%\label{parameterising point in Hy2}
%\tb{P}=\e_{12}\cosh\rho+\e_{20}\sinh\rho\cos\alpha+\e_{01}\sinh\rho\sin\alpha,
%\end{equation}
%where \(\rho\in[0,+\infty)\) and \(\alpha\in[0,2\pi)\).
%This parameterisation applies to counterclockwise points only, but then
%\(-\tb{P}\) represents a clockwise point at the same location, which is also normalised and proper.
%The point \(\tb{P}\) can also be written as 
%\begin{equation}
%\tb{P}=\e_{12}e^{\rho\e_0\wedge\tb{a}},
%\end{equation}
%where \(\tb{a}=\e_1\sin\alpha+e_2\cos\alpha\) is a normalised line perpendicular to the direction from the origin to \(\tb{P}\).
%The spinor \(S=e^{\rho\e_0\wedge\tb{a}}\) satisfies \(\e_{12}S=S^{-1}\e_{12}\).

The distance \(r\in[0,+\infty)\) between two normalised proper points \(\tb{P}\) and \(\tb{Q}\) is defined by
\begin{equation}\label{distance in Hy2}
\sinh r=\norm{\tb{P}\vee\tb{Q}}.
\end{equation}

Since \(|\tb{P}\cdot\tb{Q}|^2-\norm{\tb{P}\vee\tb{Q}}^2=1\) for the normalised proper points,
the distance \(r\) also satisfies \(\cosh r=|\tb{P}\cdot\tb{Q}|\).
%The distance measure in \Hy{2} is called hyperbolic.
%Substituting \(\tb{Q}=\e_{12}\) and \(\tb{P}\) parameterised by Equation~(\ref{parameterising point in Hy2})
%yields \(\cosh r = |\cosh\rho|\) and, therefore, \(r=\rho\) for the distance from the origin to \(\tb{P}\).
If a proper point \(\tb{P}\) is expressed via the standard coordinates \(x\) and \(y\), 
i.e.\ \(\tb{P}=\e_{12}+x\e_{20}+y\e_{01}\), then \(\norm{\tb{P}}=\sqrt{1-(x^2+y^2)}\) and the distance
from the origin to \(\tb{P}\) satisfies
\begin{equation}
\sinh{r}=\frac{\sqrt{x^2+y^2}}{\sqrt{1 - (x^2+y^2)}}, \quad\quad
\cosh{r}=\frac{1}{\sqrt{1 - (x^2+y^2)}}.
\end{equation}
If the above expressions are applied to null points, the resulting distance is infinite.
So the null points, i.e.\ points on the unit circle in \Hy{2}, are at the infinite distance from the origin.
In other words, the unit circle represents the infinity in \Hy{2}.
The improper points, which lie outside the unit circle, exist in \Hy{2} in a non-metric sense only since
they are not reachable from any proper point.
%If the two points lie on the same line passing through the origin, they can be parameterised by
%\(\tb{P}=\e_{12}e^{\rho\e_0\wedge\tb{a}}\) and \(\tb{Q}=\e_{12}e^{\eta\e_0\wedge\tb{a}}\) with a single \(\tb{a}\),
%the distance \(r\) between the points satisfies \(\cosh r=|\cosh(\rho-\eta)|\), which implies \(r=|\rho-\eta|\).

In general, the distance between proper points \(\tb{P}\) and \(\tb{Q}\) 
can be understood as the pseudo-Euclidean angle between the 1-dimensional linear subspaces in \R{3} 
corresponding to \(\J(\tb{P})\) and \(\J(\tb{Q})\), where \R{3} is viewed as a 3-dimensional Minkowski space
with the time-like direction along the \(w\)-axis.
%See Figure~\ref{visualising hyperbolic2}(a), where \(\tb{P}=\e_{12}+\tfrac{1}{2}\e_{20}-\tfrac{1}{2}\e_{01}\) 
%and \(\tb{Q}=\e_{12}-\tfrac{3}{4}\e_{20}\),  which upon the normalisation
%yields \(\cosh r=\tfrac{11/8}{\sqrt{1/2}\sqrt{7/16}}\) for the distance between the points.

\begin{figure}[t!]
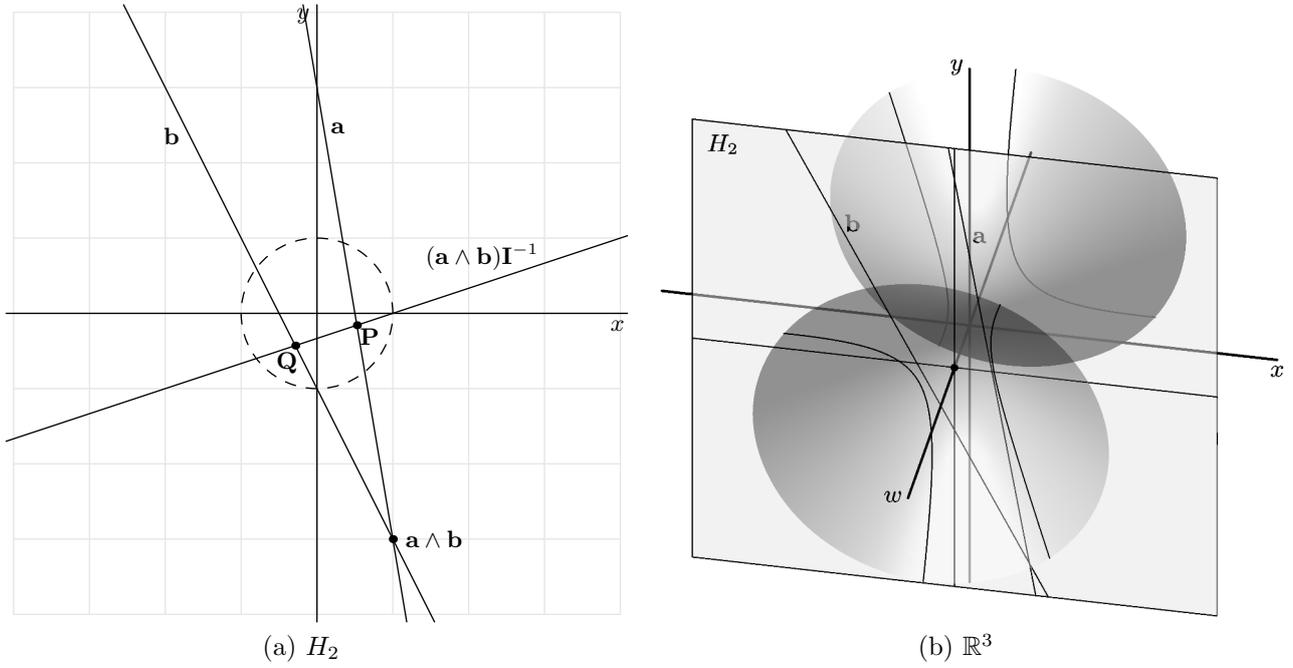

\begin{subfloatenv}{\Hy{2}}
\begin{asy}
import Figure2D;
metric=Metric(Hyperbolic);
MV a = join(Point(1,1,-3),Point(1,1/2,0) );  a/=norm(a);
MV b = join(Point(1,1,-3),Point(1,0,-1));  b/=norm(b);
MV c = wedge(a,b)/I; c/=norm(c);

MV P=wedge(a,c); 
MV Q=wedge(b,c); 

Figure f = Figure();
draw(unitcircle,dashed);

f.line(a,"$\textbf{a}$",align=(1,0),draw_orientation=false);
f.line(b,"$\textbf{b}$",draw_orientation=false);
f.line(c,"$(\textbf{a}\wedge\textbf{b})\textbf{I}^{-1}$",draw_orientation=false);

f.point(wedge(a,b), "$\textbf{a}\wedge\textbf{b}$",align=(1,0),draw_orientation=false);

f.point(P, "$\textbf{P}$",align=(0.25,-0.25),draw_orientation=false);
f.point(Q, "$\textbf{Q}$",align=(-0.25,-0.5),draw_orientation=false);

\end{asy}
\end{subfloatenv}\hspace{-12pt}%
\begin{subfloatenv}{\R{3}}
\begin{asy}
import Drawing3D;
import Figure2D;
metric=Metric(Hyperbolic);
DrawingR3 drawing = DrawingR3(4.0,0.1, camera=(10,2,7));	
drawing.target_axes();

path3 p = plane_path3((1,0,0),(1,0,0), v=(0,1,0));
drawing.Drawing3D.plane(p);
dot((1,0,0));
draw((1,-3.5,0)--(1,3.5,0));
draw((1,0,-3.5)--(1,0,3.5));

label("$H_2$", (1.5,-3,3.5));

//draw(shift((1,0,0))*rotate(90,(0,1,0))*unitcircle3,dashed);

import solids;
real r = 1;
currentlight=light(white, viewport=false,(5,-5,10));

triple f1(pair t) { return (cosh(t.x),sinh(t.x)*cos(t.y),sinh(t.x)*sin(t.y));}
triple f2(pair t) { return (-cosh(t.x),sinh(t.x)*cos(t.y),sinh(t.x)*sin(t.y));}
surface hyperboloid1=surface(f1,(0,0),(1.6,2pi),16,8,Spline);
surface hyperboloid2=surface(f2,(0,0),(1.6,2pi),16,8,Spline);
//draw(surface(hyperboloid1),lightgray+opacity(0.4),meshpen=rgb(0.6,0.6,0.6),render(compression=Low,merge=true));
//draw(surface(hyperboloid2),lightgray+opacity(0.4),meshpen=rgb(0.6,0.6,0.6),render(compression=Low,merge=true));
draw(surface(hyperboloid1),lightgray+opacity(0.4),render(compression=Low,merge=true));
draw(surface(hyperboloid2),lightgray+opacity(0.4),render(compression=Low,merge=true));

//import Figure2D;
//metric=Metric(Hyperbolic);
MV a = join(Point(1,1,-3),Point(1,1/2,0) );  a/=norm(a);
MV b = join(Point(1,1,-3),Point(1,0,-1));  b/=norm(b);

triple to3(MV P) {return (P.w,P.x,P.y);}

path3 aP = plane_path3(totriple(a), (0,0,0), scale=10);
path3 l = plane_intersection(aP,p);
draw(Label("$\textbf{a}$", 0.8,E), l);

pair[] ip = intersectionpoints(scale(tanh(1.6))*unitcircle,(1,-3)--(0,3));
MV P1 = Point(1,ip[0].x,ip[0].y);
MV Q1 = Point(1,ip[1].x,ip[1].y);
triple f(real t) { MV R = P1+(Q1-P1)*t; return to3(R/norm(R)); }
real x(real t) { return f(t).x; }
real y(real t) { return f(t).y; }
real z(real t) { return f(t).z; }
path3 aL = graph(x, y, z, 0, 1, operator ..);
draw(aL);
triple f(real t) { MV R = P1+(Q1-P1)*t; return to3(-R/norm(R)); }
real x(real t) { return f(t).x; }
real y(real t) { return f(t).y; }
real z(real t) { return f(t).z; }
path3 aL = graph(x, y, z, 0, 1, operator ..);
draw(aL);

path3 bP = plane_path3(totriple(b), (0,0,0), scale=10);
path3 l = plane_intersection(bP,p);
draw(Label("$\textbf{b}$", 0.8,E), l);

pair[] ip = intersectionpoints(scale(tanh(1.6))*unitcircle,(1,-3)--(-1,1));
MV P2 = Point(1,ip[0].x,ip[0].y);
MV Q2 = Point(1,ip[1].x,ip[1].y);
triple f(real t) { MV R = P2+(Q2-P2)*t; return to3(R/norm(R)); }
real x(real t) { return f(t).x; }
real y(real t) { return f(t).y; }
real z(real t) { return f(t).z; }
path3 bL = graph(x, y, z, 0, 1, operator ..);
draw(bL);
triple f(real t) { MV R = P2+(Q2-P2)*t; return to3(-R/norm(R)); }
real x(real t) { return f(t).x; }
real y(real t) { return f(t).y; }
real z(real t) { return f(t).z; }
path3 bL = graph(x, y, z, 0, 1, operator ..);
draw(bL);

//MV c = wedge(a,b)/I; c/=norm(c);

//MV P=wedge(a,c); 
//MV Q=wedge(b,c); 

//Figure f = Figure();
//draw(unitcircle,dashed);

//f.line(a,"$\textbf{a}$",draw_orientation=false);
//f.line(b,"$\textbf{b}$",draw_orientation=false);
//f.line(c,"$(\textbf{a}\wedge\textbf{b})\textbf{I}^{-1}$",draw_orientation=false);

//f.point(wedge(a,b), "$\textbf{a}\wedge\textbf{b}$",align=(0.25,-0.25),draw_orientation=false);

//f.point(P, "$\textbf{P}$",align=(0.25,-0.25),draw_orientation=false);
//f.point(Q, "$\textbf{Q}$",align=(0.25,-0.25),draw_orientation=false);
//P/=P.w; Q/=Q.w;

//MV ab = wedge(a,b);
//ab/=ab.w;
//dot(to3(ab));
//label("$\textbf{a}\wedge\textbf{b}$",to3(ab),(0,1,0));

//label("$\textbf{P}$",to3(P),(0,1,0));
//label("$\textbf{Q}$",to3(Q),(0,0,1));
//triple PP = to3(P/norm(P)), QQ = to3(Q/norm(Q));
//dot(to3(P)); dot(PP); dot(-PP);
//dot(to3(Q)); dot(QQ); dot(-QQ);
//draw((-PP)--PP);
//draw((-QQ)--QQ);

//triple f(real t) { MV R = P+(Q-P)*t; return to3(R/norm(R)); }
//real x(real t) { return f(t).x; }
//real y(real t) { return f(t).y; }
//real z(real t) { return f(t).z; }
//path3 PQ = graph(x, y, z, 0, 1, operator ..);
//draw(Label("$r_{min}$",(1,0.5,-1)),PQ, Arrows3(size=5));

//path3 abP = plane_path3(totriple(wedge(a,b)/I), (0,0,0), scale=10);
//path3 l = plane_intersection(abP,p);
//draw(Label("($\textbf{a}\wedge\textbf{b})\textbf{I}^{-1}$", 0.2,(2,2,0.5)), l);

\end{asy}
\end{subfloatenv}
\caption{Visualising \Hy{2}}
\label{visualising hyperbolic2}
\end{figure}

For a proper line \(\tb{a}\), its polar point  \(\tb{a}\I\) is improper.
On the other hand, if the line  \(\tb{a}\) is improper, then its polar point  \(\tb{a}\I\) is proper.
If \(\tb{a}\) is null, i.e.\ tangent to the unit circle, then \(\tb{a}\I\) is also null and is located at the point where
\(\tb{a}\) touches the unit circle.

A line intersecting \(\tb{a}\) at a proper point is perpendicular to \(\tb{a}\)
provided that it also passes through the polar point \(\tb{a}\I\)  of the line \(\tb{a}\).
The reverse is true as well, i.e.\ any line perpendicular to \(\tb{a}\) passes through \(\tb{a}\I\).
In general, the angle \(\alpha\) between two normalised proper lines \(\tb{a}\) and \(\tb{b}\) 
intersecting at a proper point is defined by
\begin{equation}
\cos\alpha=\tb{a}\cdot\tb{b}.
\end{equation}

If two proper lines intersect at an improper point, 
the distance \(r\) of the closest approach between the lines satisfies \(\cosh r=|\tb{a}\cdot\tb{b}|\)
and \(\sinh r = \norm{\tb{a}\wedge\tb{b}}\).
To find the points on the lines where this minimal separation is achieved, find
the proper line \(\tb{c}=(\tb{a}\wedge\tb{b})\I^{-1}\) whose polar point coincides with \(\tb{a}\wedge\tb{b}\). 
%the intersection point of the lines \(\tb{a}\) and \(\tb{b}\).
This line is perpendicular to both \(\tb{a}\) and \(\tb{b}\) and its intersection with \(\tb{a}\) and \(\tb{b}\)
yields the points of the closest approach, i.e.\ 
\(\cosh r=|\tb{P}\cdot\tb{Q}|\) where \(\tb{P}=\tb{a}\wedge\tb{c}/\norm{\tb{a}\wedge\tb{c}}\) 
and \(\tb{Q}=\tb{b}\wedge\tb{c}/\norm{\tb{b}\wedge\tb{c}}\).
This is illustrated in Figure~\ref{visualising hyperbolic2}(a) where 
\(\tb{a}=-\tfrac{3}{2}\e_0+3\e_1+\tfrac{1}{2}\e_2\) and
\(\tb{b}=\tfrac{1}{2}\e_0+\e_1+\tfrac{1}{2}\e_2\),
which gives \(\tb{a}\wedge\tb{b}=\e_{12}+\e_{20}-3\e_{01}\) and 
\((\tb{a}\wedge\tb{b})\I^{-1}=-\tfrac{1}{3}\e_0+\tfrac{1}{3}\e_1-\e_2\).

If two proper lines intersect at a null point, the separation between them decreases the closer one gets to the unit circle,
which represents the infinity in \Hy{2}.
In this case, one might say the lines meet at infinity.

Whether the two lines intersects at a null or an improper point, they may be considered parallel
since they do not intersect at a finite distance
(lines intersecting at improper points are also called hyperparallel).
In Euclidean and Minkowski spaces, the distance between parallel lines is constant 
(it does not depend on the position on the lines)
 and, therefore, parallel lines can be defined by means of this property.
In Hyperbolic space, this is no longer the case, since the distance between parallel lines 
depends on the position on the lines, with the minimal separation attained at unique points on the lines if they are hyperparallel.

It is possible to visualise the metric component of \Hy{2} as a  hyperboloid of two sheets defined by
\begin{equation}
w^2-x^2-y^2=1.
\end{equation}
The antipodal points on the hyperboloid, 
i.e.\ points on the sheets which lie on the same line passing through the origin,
are identified.
However, since the sheets of the hyperboloid are disconnected, it is sufficient to consider only one of them, say,
the one corresponding to \(w\ge1\).
This visualisation  is a natural extension of the metric component of \Hy{1} viewed as a hyperbola.
As in the 1-dimensional case, this visualisation has limited utility as it does not represent the distances
in a uniform fashion and does not encompass the improper points at all.
A line in \Hy{2} can be identified with a pair of parabolas 
as shown in Figure~\ref{visualising hyperbolic2}(b) for the lines \(\tb{a}\) and \(\tb{b}\).

\begin{figure}[t!]
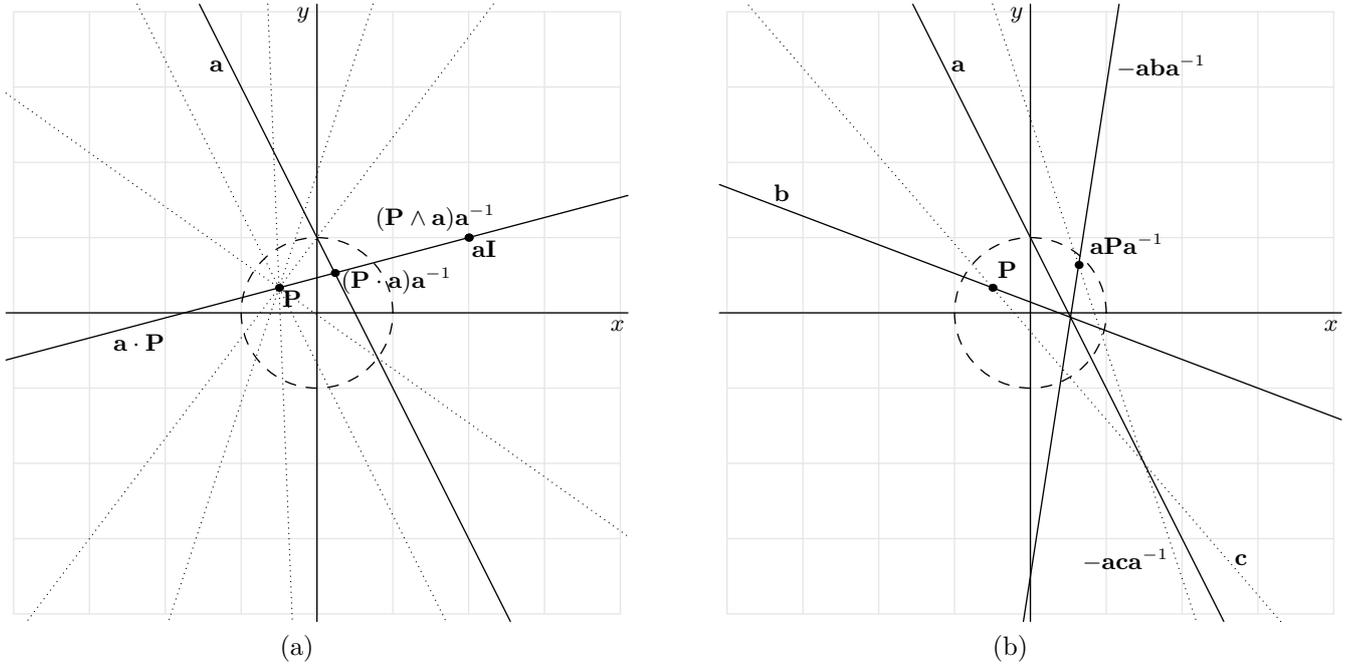

\begin{subfloatenv}{ }
\begin{asy}
import Figure2D;
Figure f = Figure();
metric = Metric(Hyperbolic);

draw(unitcircle,dashed);

MV a = join(Point(1,1/2,0),Point(1,0,1));
a=a/norm(a);
// if a^2=1
MV N1=(a-I)*dot(a,O);
MV N2=(a+I)*dot(a,O);
f.line(a,label="$\textbf{a}$",position=0.1,align=W,draw_orientation=false);
f.point(a*I,align=(0.25,-0.25),"\textbf{a}\textbf{I}",draw_orientation=false);

MV P = Point(1,-1/2,1/3);
f.point(P, "$\textbf{P}$",align=(0.25,-0.25),draw_orientation=false);
f.line(dot(a,P),"$\textbf{a}\cdot\textbf{P}$",draw_orientation=false);

f.point(dot(P,a)/a, "$(\textbf{P}\cdot\textbf{a})\textbf{a}^{-1}$",align=(0.5,-0.2),draw_orientation=false);
f.point(wedge(P,a)/a, "$(\textbf{P}\wedge\textbf{a})\textbf{a}^{-1}$",align=(-0.25,0.5),draw_orientation=false);

MV a1 = join(P,N1); a1/=norm(a1);
MV a2 = join(P,N2); a2/=norm(a2);
real alpha = pi-acos(dot(a1,a2).s);
int N = 4;
for(int i: sequence(N+1)) { MV S = exp(-1/2*alpha*i/N*P/norm(P)); f.line(S*a2/S,dotted,draw_orientation=false); }

\end{asy}
\end{subfloatenv}\hfill%
\begin{subfloatenv}{  }
\begin{asy}
import Figure2D;
Figure f = Figure();
metric = Metric(Hyperbolic);

draw(unitcircle,dashed);

MV a = join(Point(1,1/2,0),Point(1,0,1));
a=a/norm(a);
// if a^2=1
MV N1=(a-I)*dot(a,O);
MV N2=(a+I)*dot(a,O);
f.line(a,label="$\textbf{a}$",position=0.1,align=E,draw_orientation=false);
//f.point(a*I,align=(0.25,-0.25),"\textbf{a}\textbf{I}",draw_orientation=false);

MV P = Point(1,-1/2,1/3);
f.point(P, "$\textbf{P}$",align=(0.6,0.8),draw_orientation=false);

MV b = join(Point(1,-1/2,1/3),Point(1,3,-1));
b/=norm(b);
f.line(b,label="$\textbf{b}$",position=0.9,align=N,draw_orientation=false,O=topair(wedge(a,b)));
f.line(-a*b/a,label="$-\textbf{a}\textbf{b}\textbf{a}^{-1}$",position=0.1,align=E,draw_orientation=false,O=topair(wedge(a,b)));
f.point(a*P/a,"$\textbf{a}\textbf{P}\textbf{a}^{-1}$",align=(0.9,0.9),draw_orientation=false);

MV R = Point(1,1.5,-2);
MV c = join(P,R);
c/=norm(c);
f.line(c,currentpen+dotted,label="$\textbf{c}$",position=0.1,align=E,draw_orientation=false,O=topair(R));
f.line(-a*c/a,currentpen+dotted,label="$-\textbf{a}\textbf{c}\textbf{a}^{-1}$",position=0.9,align=W,draw_orientation=false,O=topair(R));

\end{asy}
\end{subfloatenv}
\caption{Basic properties of \Hy{2}}
\label{orientation and polar points in H2}
\end{figure}

The polar point \(\tb{a}\I=\tfrac{1}{2}\e_{12}+ \e_{20}+\tfrac{1}{2}\e_{01}\) 
of the proper line \(\tb{a}=-\tfrac{1}{2}\e_0+ \e_1+\tfrac{1}{2}\e_2\) is shown in
Figure~\ref{orientation and polar points in H2}(a).
The null points where \(\tb{a}\) intersects the unit circle are given by \(\tb{N}_{\pm}=(\tb{a}\pm\I)(\tb{a}\cdot\e_{12})\)
if the proper line \(\tb{a}\) is normalised.
The lines passing through \(\tb{a}\I\) and the points \(\tb{N}_+\) and \(\tb{N}_-\) are tangent to the unit circle, i.e.\ 
\((\tb{a}\I)\vee\tb{N}_+\)  and \((\tb{a}\I)\vee\tb{N}_-\) are null.
For any two proper points \(\tb{P}\) and \(\tb{Q}\), the commutator \(\tb{P}\times\tb{Q}\) gives the
polar point of the line \(\tb{P}\vee\tb{Q}\).

The line passing through the proper point \(\tb{P}=\e_{12}-\tfrac{1}{2}\e_{20}+\tfrac{1}{3}\e_{01}\) 
and perpendicular to \(\tb{a}\) is given by \(\tb{a}\cdot\tb{P}\)
as usual.
Note that \(\tb{a}\cdot\tb{P}\) also passes through the polar point of \(\tb{a}\).
There are infinite number of lines in \Hy{2} which pass through \(\tb{P}\) and are parallel to \(\tb{a}\);
some of them are shown in Figure~\ref{orientation and polar points in H2}(a) with dotted lines.
The projection of \(\tb{P}\) on \(\tb{a}\) is defined in the usual way and is located at the intersection
of \(\tb{a}\) and \(\tb{a}\cdot\tb{P}\).
The rejection is located at the polar point of \(\tb{a}\) as expected.

The distance \(r\) between a proper line \(\tb{a}\) and a proper point \(\tb{P}\) satisfies
\(\sinh r=|\tb{a}\vee\tb{P}|\) provided that the line and the point are normalised.
The same distance also satisfies \(\cosh r=\norm{\tb{a}\cdot\tb{P}}\) since 
\(\norm{\tb{a}\cdot\tb{P}}^2-  |\tb{a}\vee\tb{P}|^2=1\) in \Hy{2}.
For instance, \(\sinh r=\tfrac{5/6}{\sqrt{23}/6}\) for \(\tb{a}\) and \(\tb{P}\) shown in 
Figure~\ref{orientation and polar points in H2}(a).
For any proper line \(\tb{a}=d\e_0+a\e_1+b\e_2\), which is normalised,
the distance from the origin to the line satisfies \(\sinh r=|d|\).

The area \({\cal S}\) of a right triangle in \Hy{2}
defined by the proper points \(\tb{P}\), \(\tb{Q}\), \(\tb{R}\) satisfies
the same formula as in elliptic space \El{2}. Namely,
\begin{equation}
\sin{\cal S}=\frac{|\tb{P}\vee\tb{Q}\vee\tb{R}|}{1+|\tb{Q}\cdot\tb{R}|}
\end{equation}
where the points are assumed to be normalised and \((\tb{P}\vee\tb{Q})\cdot(\tb{P}\vee\tb{R})=0\),
i.e.\ the angle at \(\tb{P}\) equals~\(\tfrac{\pi}{2}\).
However, since the metric is different, the numeric value of the area is different from that in \El{2}
for a given triangle.
The above expression is not applicable if either \(\tb{Q}\) or \(\tb{R}\) is null 
(\(\tb{P}\) cannot be null, otherwise the angle at \(\tb{P}\) would be zero rather than \(\tfrac{\pi}{2}\)).
By treating the normalisation explicitly one obtains an equivalent expression:
\(\sin{\cal S}=|\tb{P}\vee\tb{Q}\vee\tb{R}|/(\norm{\tb{Q}}\norm{\tb{R}}+  |\tb{Q}\cdot\tb{R}|)\),
where \(\tb{Q}\) and \(\tb{R}\) need not be normalised and, therefore, can be null, in which case
\(\sin{\cal S}=|\tb{P}\vee\tb{Q}\vee\tb{R}|/|\tb{Q}\cdot\tb{R}|\) where \(\tb{Q}\) and  \(\tb{R}\)
are not normalised but \(\tb{P}\) is.

The area \(\cal S\) of a general hyperbolic triangle is given by 
\begin{equation}
{\cal S}=-\alpha-\beta+\gamma,
\end{equation}
which is also equal to  \({\cal S}=\pi-(\alpha+\beta+\pi-\gamma)\),
where 
\(\alpha=\arccos(\frac{\tb{r}\cdot\tb{q}}{\norm{\tb{r}}\norm{\tb{q}}})\),
\(\beta=\arccos(\frac{\tb{r}\cdot\tb{p}}{\norm{\tb{r}}\norm{\tb{p}}})\),
\(\gamma=\arccos(\frac{\tb{q}\cdot\tb{p}}{\norm{\tb{q}}\norm{\tb{p}}})\),
and \(\tb{r}=\tb{P}\vee\tb{Q}\), \(\tb{q}=\tb{P}\vee\tb{R}\), \(\tb{p}=\tb{R}\vee\tb{Q}\).
Note also that even though hyperbolic space \Hy{2} is unbound and has an infinite area,
 there is an upper bound 
on the area a triangle in \Hy{2} can assume. 
It is equal to \(\pi\), which is the area of any triangle whose corners lie at infinity
(\(\alpha,\beta=0\) and \(\gamma=\pi\) in this case).
The upper bound on the area of the right triangle is \(\tfrac{\pi}{2}\).

%import Figure2D;
%metric = Metric(Hyperbolic);
%MV a,b,c, P,Q,R;
%pair random_pair_in_square() {return (1-2*unitrand(),1-2*unitrand()); }
%bool outside_unitcircle(pair p) { return length(p)>1; }
%MV random_point() { pair p; do {p=random_pair_in_square();} while(outside_unitcircle(p)); return Point(1,p.x,p.y);}
%struct T {MV P,Q,R; };
%T t = new T;
%real A, maxA;
%for(int i:sequence(10000)) {
%  P = random_point(); Q = random_point(); R = random_point();
%  a = join(P,Q); b = join(P,R); c = join(Q,R);
%  a/=norm(a); b/=norm(b); c/=norm(c);
%  real alpha = acos(dot(a,b).s);
%  real beta = acos(dot(b,c).s);
%  real gamma = acos(dot(a,c).s); 
%  A = gamma - alpha - beta;
%  //A = asin(abs(join(P,Q,R).s)/(1+abs(dot(Q,R).s)));
%  if(A>maxA) { maxA = A; write(maxA); t.P = P; t.Q = Q; t.R = R; }
%}
%Figure f = Figure();
%draw(unitcircle,dashed);
%f.point(t.P,'P',draw_orientation=false);
%f.point(t.Q,'Q',draw_orientation=false);
%f.point(t.R,'R',draw_orientation=false);
%a = join(t.P,t.Q); b = join(t.P,t.R); c = join(t.Q,t.R);
%a/=norm(a); b/=norm(b); c/=norm(c);
%f.line(a,'a',draw_orientation=false,draw_bottom_up_orientation=true,O=topair(t.P));
%f.line(a,'',draw_orientation=false,draw_bottom_up_orientation=true,O=topair(t.Q));
%f.line(b,'b',draw_orientation=false,draw_bottom_up_orientation=true,O=topair(t.P));
%f.line(b,'',draw_orientation=false,draw_bottom_up_orientation=true,O=topair(t.R));
%f.line(c,'c',draw_orientation=false,draw_bottom_up_orientation=true,O=topair(t.Q));
%f.line(c,'',draw_orientation=false,draw_bottom_up_orientation=true,O=topair(t.R));

\begin{figure}[t!]
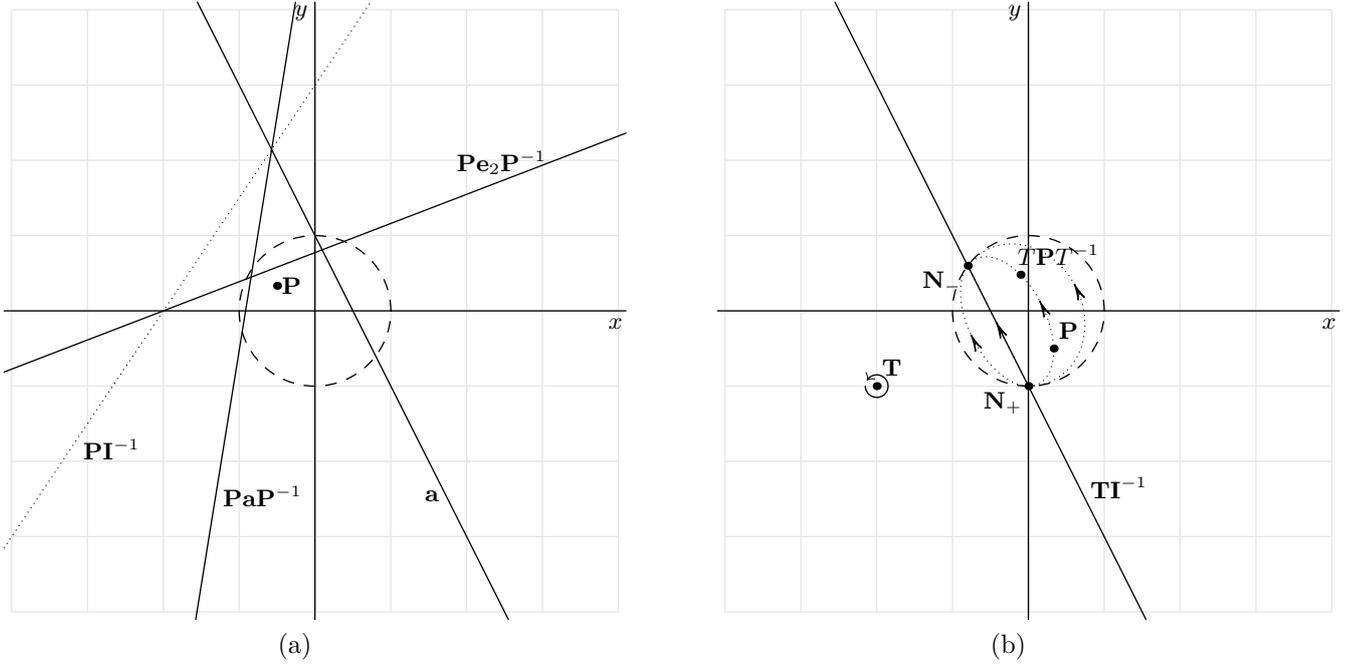

\begin{subfloatenv}{ }
\begin{asy}
import Figure2D;
Figure f = Figure();
metric = Metric(Hyperbolic);

draw(unitcircle,dashed);

MV a = join(Point(1,1/2,0),Point(1,0,1));
a=a/norm(a);
// if a^2=1
MV N1=(a-I)*dot(a,O);
MV N2=(a+I)*dot(a,O);
f.line(a,label="$\textbf{a}$",position=0.8,align=W,draw_orientation=false);

MV P = Point(1,-1/2,1/3);
f.point(P, "$\textbf{P}$",align=(0.4,0),draw_orientation=false);

f.line(P*a/P,label="$\textbf{P}\textbf{a}\textbf{P}^{-1}$",position=0.2,align=E,draw_orientation=false);
//f.line(P*e_1/P,label="$\textbf{P}\textbf{e}_1\textbf{P}^{-1}$",position=0.9,align=E,draw_orientation=false);
f.line(P*e_2/P,label="$\textbf{P}\textbf{e}_2\textbf{P}^{-1}$",position=0.2,align=N,draw_orientation=false);
f.line(P/I,currentpen+dotted,"$\textbf{P}\textbf{I}^{-1}$",draw_orientation=false);

\end{asy}
\end{subfloatenv}\hfill%
\begin{subfloatenv}{  }
\begin{asy}
import Figure2D;
Figure f = Figure();
metric = Metric(Hyperbolic);

draw(unitcircle,dashed);

MV T = Point(1,-2,-1);
T/=norm(T);
MV a = join(O,T);
a=a/norm(a);
// if a^2=1
MV N=(a+I)*dot(a,O);

f.point(T, "$\textbf{T}$",align=(0.6,0.8),draw_orientation=true);
f.line(T/I,"$\textbf{T}\textbf{I}^{-1}$",draw_orientation=false);

MV P = Point(1,1/3,-1/2);
f.point(P, "$\textbf{P}$",align=(0.6,0.8),draw_orientation=false);

MV N1=(T+1)*cross(T,O);
MV N2=(T-1)*cross(T,O);
f.point(N1, "$\textbf{N}_+$",align=(-0.5,-0.5),draw_orientation=false);
f.point(N2, "$\textbf{N}_-$",align=(-0.5,-0.5),draw_orientation=false);

real lambda=1;
MV S = exp(-1/2*lambda*T);
f.point(S*P/S, "$T\textbf{P}T^{-1}$",align=(0.6,0.8),draw_orientation=false);
real arrow_size=3;
real x(real lambda) {MV S = exp(-1/2*lambda*T); return topair(S*P/S).x; }
real y(real lambda) {MV S = exp(-1/2*lambda*T); return topair(S*P/S).y; }
draw(graph(x,y,-5,5),p=dotted,MidArrow(HookHead,size=arrow_size));

MV a = -dot(T/I,P);
a/=norm(a);
P/=norm(P);
real lam =  asinh(abs(join(T/I,P).s));
MV S = exp(-1/2*lam*a*I);
MV P2 = S*P/S;
real x(real lambda) {MV S = exp(-1/2*lambda*T); return topair(S*P2/S).x; }
real y(real lambda) {MV S = exp(-1/2*lambda*T); return topair(S*P2/S).y; }
draw(graph(x,y,-5,5),p=dotted,MidArrow(HookHead,size=arrow_size));

MV S = exp(-1/2*(-lam)*a*I);
MV P2 = S*P/S;
real x(real lambda) {MV S = exp(-1/2*lambda*T); return topair(S*P2/S).x; }
real y(real lambda) {MV S = exp(-1/2*lambda*T); return topair(S*P2/S).y; }
draw(graph(x,y,-5,5),p=dotted,MidArrow(HookHead,size=arrow_size));

MV S = exp(-1/2*(-2*lam)*a*I);
MV P2 = S*P/S;
real x(real lambda) {MV S = exp(-1/2*lambda*T); return topair(S*P2/S).x; }
real y(real lambda) {MV S = exp(-1/2*lambda*T); return topair(S*P2/S).y; }
draw(graph(x,y,-5,5),p=dotted,MidArrow(HookHead,size=arrow_size));

\end{asy}
\end{subfloatenv}
\caption{Reflection and translation in \Hy{2}}
\label{reflection in a point H2}
\end{figure}

Reflection in a proper line is defined as usual and is illustrated in Figure~\ref{orientation and polar points in H2}(b).
Note that \(\tb{a}\tb{P}\tb{a}^{-1}\) lies on the line \(\tb{a}\cdot\tb{P}\), which is perpendicular to \(\tb{a}\).
The reflection of \(\tb{b}\), which intersects \(\tb{a}\) at a proper point, preserves the angle between the lines,
whereas the reflection of \(\tb{c}\), which intersects \(\tb{a}\) at an improper point, preserves the minimal distance
between the lines.

The reflection in a proper point is illustrated in Figure~\ref{reflection in a point H2}(a) where
\(\tb{P} = \e_{12}-\tfrac{1}{2}\e_{20}+\tfrac{1}{3}\e_{01}\).
The reflection of a line in a point preserves the distance from the point to the line.
The intersection point of the original line and the reflected line is improper and lies on the line \(\tb{P}\I^{-1}\),
which is shown in Figure~\ref{reflection in a point H2}(a) with the dotted line.

The spin group of \Hy{2} is defined in the usual ways and consists of spinors of the form \(S=e^B\) where \(B\) is a bivector.
The action of \(SA_kS^{-1}\) of the spinor \(S\) on the geometric object dually represented by \(A_k\)
depends on whether the point dually represented by \(B\) is proper, null, or improper.

The action of  \(T=e^{-\tfrac{1}{2}\lambda\tb{T}}\), where the point \(\tb{T}\) is normalised and improper (\(\tb{T}^2=1\)), 
is equivalent to a translation by \(\lambda\) along the proper line \(\tb{T}\I^{-1}\).
An example is shown in Figure~\ref{reflection in a point H2}(b) where 
\(\tb{T}=\tfrac{1}{2}\e_{12}-\e_{20}-\tfrac{1}{2}\e_{01}\), 
\(\tb{P}=\e_{12}+\tfrac{1}{3}\e_{20}-\tfrac{1}{2}\e_{01}\), and \(\lambda=1\).
The trajectory of \(\tb{P}\) under the action of \(T\) with different values of \(\lambda\in\R{}\) is shown with
a dotted line passing through \(\tb{P}\); the trajectories for several other points are also shown.
Even though it appears that all such trajectories depart from the point \(\tb{N}_+=(\tb{T}+1)(\tb{T}\times\e_{12})\)
and terminate at \(\tb{N}_-=(\tb{T}-1)(\tb{T}\times\e_{12})\), where \(\tb{T}\) is assumed to be normalised,
the proper points located on a given trajectory are at the same distance from the line \(\tb{T}\I^{-1}\).
%So the translation induced by the spinor \(T\) can be thought of as a translation along the line \(\tb{T}\I^{-1}\).
The null points \(\tb{N}_\pm\) remain fixed under the action of \(T\) but their weight changes;
all other null points move along the unit circle towards \(\tb{N}_-\).
Proper points located on the line \(\tb{T}\I^{-1}\) remain on this line under the translation and the line \(\tb{T}\I^{-1}\) itself is invariant under the translation.

\begin{figure}[t!]
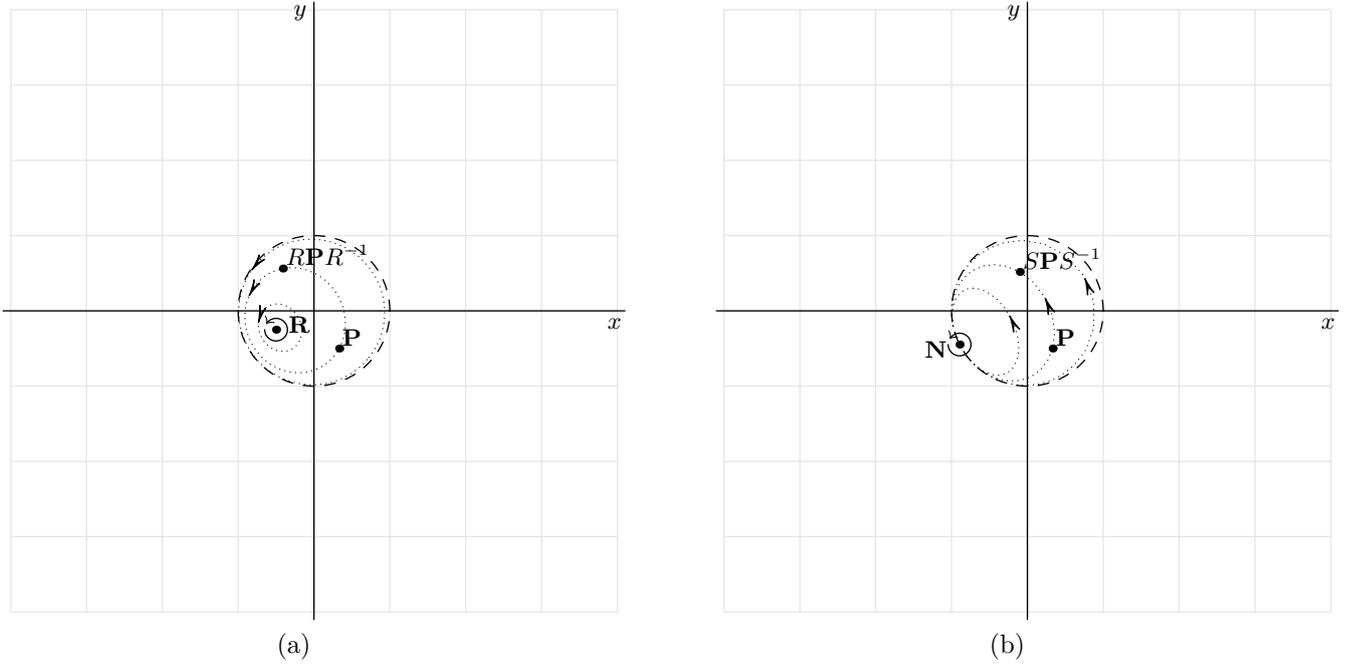

\begin{subfloatenv}{ }
\begin{asy}
import Figure2D;
Figure f = Figure();
metric = Metric(Hyperbolic);

draw(unitcircle,dashed);

MV R = Point(1,-1/2,-1/4);
R/=norm(R);

f.point(R, "$\textbf{R}$",align=(1,0.25),draw_orientation=true);

MV P = Point(1,1/3,-1/2);
f.point(P, "$\textbf{P}$",align=(0.25,0.25),draw_orientation=false);

real lambda=pi/2;
MV S = exp(-1/2*lambda*R);
f.point(S*P/S, "$R\textbf{P}R^{-1}$",align=(0.25,0.25),draw_orientation=false);
real arrow_size=3;
real x(real lambda) {MV S = exp(-1/2*lambda*R); return topair(S*P/S).x; }
real y(real lambda) {MV S = exp(-1/2*lambda*R); return topair(S*P/S).y; }
draw(graph(x,y,0,2*pi),p=dotted,MidArrow(HookHead,size=arrow_size));

P/=norm(P);
real lam =  asinh(norm(join(P,R)));
MV S = exp(-1/2*lam*join(P,R)*I);
MV P2 = S*P/S;
real x(real lambda) {MV S = exp(-1/2*lambda*R); return topair(S*P2/S).x; }
real y(real lambda) {MV S = exp(-1/2*lambda*R); return topair(S*P2/S).y; }
draw(graph(x,y,0,2*pi),p=dotted,MidArrow(HookHead,size=arrow_size));

MV S = exp(-1/2*(-lam/2)*join(P,R)*I);
MV P2 = S*P/S;
real x(real lambda) {MV S = exp(-1/2*lambda*R); return topair(S*P2/S).x; }
real y(real lambda) {MV S = exp(-1/2*lambda*R); return topair(S*P2/S).y; }
draw(graph(x,y,0,2*pi),p=dotted,MidArrow(HookHead,size=arrow_size));

\end{asy}
\end{subfloatenv}\hfill%
\begin{subfloatenv}{  }
\begin{asy}
import Figure2D;
Figure f = Figure();
metric = Metric(Hyperbolic);

draw(unitcircle,dashed);

MV R = Point(1,-1/2,-1/4);
MV a = join(O,R);
a/=norm(a);
MV N = (a+I)*dot(a,O);

f.point(N, "$\textbf{N}$",align=(-1,-0.25),draw_orientation=true);

MV P = Point(1,1/3,-1/2);
f.point(P, "$\textbf{P}$",align=(0.25,0.25),draw_orientation=false);

real lambda=1;
MV S = exp(-1/2*lambda*N);
f.point(S*P/S, "$S\textbf{P}S^{-1}$",align=(0.25,0.25),draw_orientation=false);
real arrow_size=3;
real x(real lambda) {MV S = exp(-1/2*lambda*N); return topair(S*P/S).x; }
real y(real lambda) {MV S = exp(-1/2*lambda*N); return topair(S*P/S).y; }
draw(graph(x,y,-10,10),p=dotted,MidArrow(HookHead,size=arrow_size));

real lam = 1;
MV S = exp(-1/2*lam*join(P,N)*I);
MV P2 = S*P/S;
real x(real lambda) {MV S = exp(-1/2*lambda*N); return topair(S*P2/S).x; }
real y(real lambda) {MV S = exp(-1/2*lambda*N); return topair(S*P2/S).y; }
draw(graph(x,y,-10,10),p=dotted,MidArrow(HookHead,size=arrow_size));

MV S = exp(-1/2*(-lam/2)*join(P,N)*I);
MV P2 = S*P/S;
real x(real lambda) {MV S = exp(-1/2*lambda*N); return topair(S*P2/S).x; }
real y(real lambda) {MV S = exp(-1/2*lambda*N); return topair(S*P2/S).y; }
draw(graph(x,y,-10,10),p=dotted,MidArrow(HookHead,size=arrow_size));

\end{asy}
\end{subfloatenv}
\caption{Rotation and null translation in \Hy{2}}
\label{rotation in H2}
\end{figure}

The action of  \(R=e^{-\tfrac{1}{2}\alpha\tb{R}}\), 
where the point \(\tb{R}\) is normalised and proper (\(\tb{R}^2=-1\)), 
is equivalent to a rotation around \(\tb{R}\) by the angle \(\alpha\).
An example is given in Figure~\ref{rotation in H2}(a),
 where 
\(\tb{P}=\e_{12}+\tfrac{1}{3}\e_{20}-\tfrac{1}{2}\e_{01}\),
\(\tb{R}=\tfrac{1}{\sqrt{11}}(4\e_{12}-2\e_{20}-1\e_{01})\), and \(\alpha=\tfrac{\pi}{2}\).
The trajectory of \(\tb{P}\) under the action of \(R\) with different values of \(\alpha\) 
is shown with a dotted line passing through \(\tb{P}\);
the trajectories for two more points are shown as well.
All points along a given trajectory are at the same distance from the point \(\tb{R}\).
So the trajectories represent circles in \Hy{2} with the centre at \(\tb{R}\).
The centre of rotation \(\tb{R}\) is invariant under the action of the spinor \(R\).
The angle \(\alpha=2\pi\) is required to complete one full rotation around a hyperbolic circle.

Any null point in \Hy{2} can be written as  \(\theta\tb{N}\), where \(\theta\in\R{}\),
 \(\tb{N}=\e_{12}+\e_0\wedge\tb{a}\), and \(\tb{a}\) is a normalised proper line passing through the origin and \(\tb{N}\).
The action of  \(S=e^{-\tfrac{1}{2}\theta\tb{N}}\), where \(\tb{N}\) is expressed as above, 
may be called a null translation by \(\theta\) anchored to the null point \(\tb{N}\).
Since \(\tb{N}^2=0\), the action of the spinor \(S\) on a proper point \(\tb{P}\) simplifies to 
\begin{equation}
e^{-\tfrac{1}{2}\theta\tb{N}}\tb{P}
e^{\tfrac{1}{2}\theta\tb{N}}
=\tb{P}+\theta \tb{P}\times\tb{N} -\tfrac{1}{4}\theta^2 \tb{N}\tb{P}\tb{N}.
\end{equation}
It is illustrated in Figure~\ref{rotation in H2}(b) where \(\tb{N}=\e_{12}+\e_0\wedge\tb{a}\),
\(\tb{a}=\tfrac{1}{\sqrt{5}}(-\e_1+2\e_2)\), and \(\theta=1\).
The trajectories generated by \(S\) with different values of \(\theta\in\R{}\)
for three points including \(\tb{P}\) are shown  with the dotted lines.
None of the proper points or proper lines is invariant under the null translation,
but the null point, which the null translation is anchored to, is invariant.

\section{Hyperbolic space \Hy{3}}

For a plane \(\tb{a}=d\e_0+a\e_1+b\e_2+c\e_3\), 
a line \(\mb{\Lambda}=p_{10}\e_{10}+p_{20}\e_{20}+p_{30}\e_{30}+p_{23}\e_{23}+p_{31}\e_{31}+p_{12}\e_{12}\),
where \(p_{10}p_{23}+p_{20}p_{31}+ p_{30}p_{12}=0\), and 
a point \(\tb{P}=w\e_{123}+x\e_{320}+y\e_{130}+z\e_{210}\), 
I compute  \(\tb{a}^2=-d^2+a^2+b^2+c^2\), 
\(\mb{\Lambda}^2=p_{10}^2+p_{20}^2+p_{30}^2-p_{23}^2-p_{31}^2-p_{12}^2\),
and \(\tb{P}^2=-w^2+x^2+y^2+z^2\), which gives
\(\norm{\tb{a}}=\sqrt{|\tb{a}^2|}\), \(\norm{\mb{\Lambda}}=\sqrt{|\mb{\Lambda}^2|}\),
and \(\norm{\tb{P}}=\sqrt{|\tb{P}^2|}\).
A plane \(\tb{a}\) is called proper if \(\tb{a}^2>0\), null if \(\tb{a}^2=0\), and improper otherwise.
A line \(\mb{\Lambda}\) is proper if \(\mb{\Lambda}^2<0\), null if \(\mb{\Lambda}^2=0\),
and improper otherwise.
A point \(\tb{P}\) is proper if \(\tb{P}^2<0\), null if \(\tb{P}^2=0\), and improper otherwise.
The null points comprise the unit\footnote{Once again, this terminology is used in a non-metric sense.}  sphere centred on the origin of \Hy{3};
 the proper points are located inside this sphere and improper points are outside.
Proper lines and planes pass through some proper points and
null lines and planes are tangent to the unit sphere.
Hence, the metric content of \Hy{3} is confined to the inside of the unit sphere;
all other points exist in \Hy{3} in a non-metric sense only.

\begin{figure}[t]
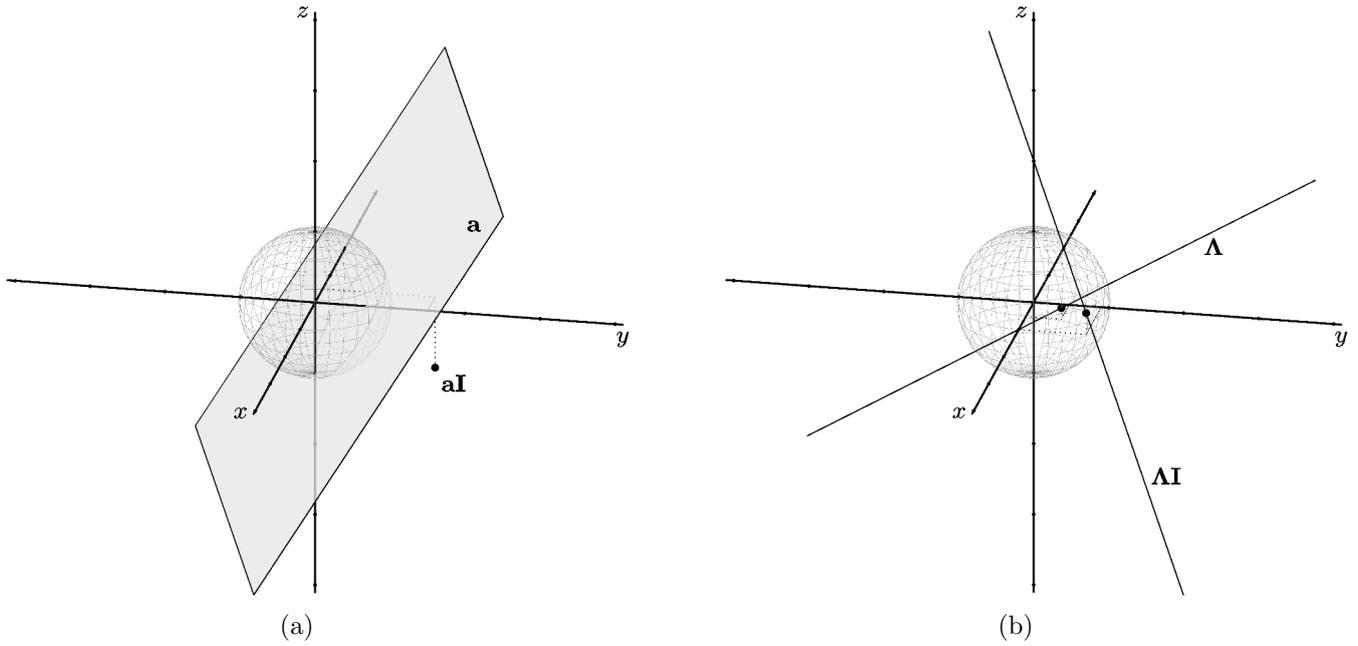
\hspace{-1cm}
\begin{subfloatenv}{ }
\begin{asy}
import Figure3D;
Figure f = Figure();
metric=Metric(Hyperbolic);

import solids;
//triple f(pair t) { return (cos(t.x),sin(t.x)*cos(t.y),sin(t.x)*sin(t.y));}
triple f(pair t) { return (sin(t.x)*cos(t.y),sin(t.x)*sin(t.y),cos(t.x));}
surface sphere=surface(f,(0,0),(pi,2pi),16,32,Spline);
draw(surface(sphere),nullpen,meshpen=rgb(0.6,0.6,0.6),render(compression=Low,merge=true));

var a = Plane(1,1/2,-3/2,1);
f.plane(a, "$\textbf{a}$", align=(1,-3,-1),draw_orientation=false);

triple c = totriple(a.centre());
real r = sqrt(1-length(c)^2);
draw(circle(c,r,c),rgb(0.6,0.6,0.6)+0.15);

//var L = wedge(e_0,a);
//f.line_at_infinity(L, centre=(0,-1,1.5), label="$\textbf{e}_0\wedge\textbf{a}$", align=(0,-1,0));

f.point(a*I, "$\textbf{a}\textbf{I}$", draw_orientation=false, align=(0,1,-1));

\end{asy}
\end{subfloatenv}\hfill%
\begin{subfloatenv}{ }
\begin{asy}
import Figure3D;
Figure f = Figure();
metric=Metric(Hyperbolic);

import solids;
triple f(pair t) { return (sin(t.x)*cos(t.y),sin(t.x)*sin(t.y),cos(t.x));}
surface sphere=surface(f,(0,0),(pi,2pi),16,32,Spline);
draw(surface(sphere),nullpen,meshpen=rgb(0.6,0.6,0.6),render(compression=Low,merge=true));

var P = Point(1,1,0,0);
var Q = Point(1,0,1,1/3);
var L = join(P,Q);
write(L);

f.line(L, "$\boldsymbol{\Lambda}$", align=(0,0,-1),draw_orientation=false);

MV C = toline(L).centre();
f.point(C,"",draw_orientation=false,draw_helper_lines=true);

f.line(L*I, "$\boldsymbol{\Lambda}\textbf{I}$", position=0.8,align=(0,1,0.5),draw_orientation=false);

MV CI = toline(L*I).centre();
f.point(CI,"",draw_orientation=false,draw_helper_lines=true);

\end{asy}
\end{subfloatenv}
\caption{Basic properties of points, lines, and planes  in \Hy{3}}
\label{basic  H3}
\end{figure}

The distance \(r\) between two normalised proper points \(\tb{P}\) and \(\tb{Q}\) is defined by
\begin{equation}
\sinh r = \norm{\tb{P}\vee\tb{Q}}
\end{equation}
and is also given by \(\cosh r = |\tb{P}\cdot\tb{Q}|\).
The null points are at the infinite distance from any proper point.
The distance to or between improper points is not defined.
If a plane \(\tb{a}\) is proper then its polar point \(\tb{a}\I\) is improper, and vice versa
(see Figure~\ref{basic  H3}(a) where the unit sphere is shown as a wire frame).
If \(\tb{a}\) is null then its polar point is also null and it is located at the point where \(\tb{a}\)
touches the unit sphere.
If a line \(\mb{\Lambda}\) is proper then \(\mb{\Lambda}\I\) is improper, and vice versa
(see Figure~\ref{basic  H3}(b)).
The null points where a proper line \(\mb{\Lambda}\) intersects the unit sphere are given by 
\(\tb{N}_\pm=(\mb{\Lambda}\pm\I)(\mb{\Lambda}\cdot\e_{123})\) provided that \(\mb{\Lambda}\) is normalised.
If \(\mb{\Lambda}\) is null then \(\mb{\Lambda}\I\) is also null and it touches the unit sphere at
the same point as \(\mb{\Lambda}\).
The null point where a null line  \(\mb{\Lambda}\) touches the unit sphere is given by \(\mb{\Lambda}(\mb{\Lambda}\cdot\e_{123})\).
The angle \(\alpha\) between  two normalised proper planes
\(\tb{a}\) and \(\tb{b}\) is defined by 
\begin{equation}
\cos\alpha=\tb{a}\cdot\tb{b},
\end{equation}
provided that the planes intersect along a proper line.
Otherwise the angle between planes is undefined.
Any plane perpendicular to a proper plane \(\tb{a}\) passes through \(\tb{a}\I\).
Similarly, any plane perpendicular to a proper line \(\mb{\Lambda}\) passes through \(\mb{\Lambda}\I\).

The distance \(r\) from a normalised proper point \(\tb{P}\) to a normalised proper plane \(\tb{a}\) 
satisfies \(\sinh r=|\tb{a}\vee\tb{P}|\) and \(\cosh r = \norm{\tb{a}\cdot\tb{P}}\).
It equals the distance along the line \(\tb{a}\cdot\tb{P}\) from \(\tb{P}\) to the point where it intersects \(\tb{a}\).
The distance from a normalised proper point \(\tb{P}\) 
to a normalised proper line \(\mb{\Lambda}\) satisfies \(\sinh r=\norm{\mb{\Lambda}\vee\tb{P}}\)
and \(\cosh r = \norm{\mb{\Lambda}\cdot\tb{P}}\), 
where \(\mb{\Lambda}\vee\tb{P}\) is a plane passing through \(\mb{\Lambda}\) and \(\tb{P}\),
and \(\mb{\Lambda}\cdot\tb{P}\) is a plane perpendicular to \(\mb{\Lambda}\) and passing through \(\tb{P}\).
The angle \(\alpha\) between a line \(\mb{\Lambda}\) and a plane \(\tb{a}\) 
intersecting at a proper point
satisfies \(\cos\alpha=\norm{\tb{a}\cdot\mb{\Lambda}}\) and \(\sin\alpha=\norm{\tb{a}\wedge\mb{\Lambda}}\),
where \(\tb{a}\cdot\mb{\Lambda}\) is a plane passing through \(\mb{\Lambda}\) and perpendicular to \(\tb{a}\),
and \(\tb{a}\wedge\mb{\Lambda}\) is a point where \(\mb{\Lambda}\) intersects \(\tb{a}\)
(the angle \(\alpha\) does not take into account the orientation of \(\mb{\Lambda}\) and \(\tb{a}\)).
The angle \(\alpha\) between two normalised proper lines \(\mb{\Lambda}\) and \(\mb{\Theta}\) 
intersecting at a proper point is given by
\(\mb{\Lambda}\cdot\mb{\Theta}=-\cos\alpha\).
If \(\mb{\Lambda}\) and \(\mb{\Theta}\) intersect at an improper point,
then the distance \(r\) between them satisfies \(|\mb{\Lambda}\cdot\mb{\Theta}|=\cosh r\).

Any non-simple bivector \(\mb{\Lambda}\) in \Hy{3} can be decomposed into a sum of two complementary lines,
\(\mb{\Lambda}=\mb{\Lambda}_1+\mb{\Lambda}_2\),
called the axes of \(\mb{\Lambda}\), one of which is proper and the other is improper.
This follows from the fact that \((\mb{\Lambda}\wedge\mb{\Lambda})^2<0\) for any non-simple bivector \mb{\Lambda}.
The axes  are given by
\begin{equation}
\mb{\Lambda}_1=\mb{\Lambda}/(1+\tfrac{1}{2}(\mb{\Lambda}\wedge\mb{\Lambda})/\mb{\Lambda}_1^2)\quad
\text{and}\quad
\mb{\Lambda}_2=\mb{\Lambda}/(1+\tfrac{1}{2}(\mb{\Lambda}\wedge\mb{\Lambda})/\mb{\Lambda}_2^2),
\end{equation}
where
\begin{equation}
\mb{\Lambda}_{1,2}^2=\tfrac{1}{2}\left(\mb{\Lambda}\cdot\mb{\Lambda}
\pm\sqrt{(\mb{\Lambda}\cdot\mb{\Lambda})^2-(\mb{\Lambda}\wedge\mb{\Lambda})^2}\right)
\end{equation}
or
\begin{equation}
\mb{\Lambda}_{1,2}^2=\tfrac{1}{2}\left(\mb{\Lambda}\cdot\mb{\Lambda}
\pm\sqrt{(\mb{\Lambda}\cdot\mb{\Lambda})^2+(\mb{\Lambda}\vee\mb{\Lambda})^2}\right),
\label{L1^2 and L2^2}
\end{equation}
since \(\I^2=-1\).
For definiteness, I will assume \(\mb{\Lambda}_1\) uses the minus sign in (\ref{L1^2 and L2^2}) and \(\mb{\Lambda}_2\) uses the plus sign,
so that \(\mb{\Lambda}_1\) is proper and  \(\mb{\Lambda}_2\) is improper.
I get \(\mb{\Lambda}_2/\norm{\mb{\Lambda}_2}=\mb{\Lambda}_1\I/\norm{\mb{\Lambda}_1}\)
if \(\mb{\Lambda}\vee\mb{\Lambda}<0\) and 
\(\mb{\Lambda}_2/\norm{\mb{\Lambda}_2}=-\mb{\Lambda}_1\I/\norm{\mb{\Lambda}_1}\)
if \(\mb{\Lambda}\vee\mb{\Lambda}>0\).
If \(\mb{\Lambda}\cdot\mb{\Lambda}=0\) then the above formulas simplify to  \(\mb{\Lambda}_{1,2}=\tfrac{1}{2}(1\pm\I)\mb{\Lambda}\) (note that \(1+\I\) and \(1-\I\) are invertible in \Hy{3}).

\begin{figure}[t]
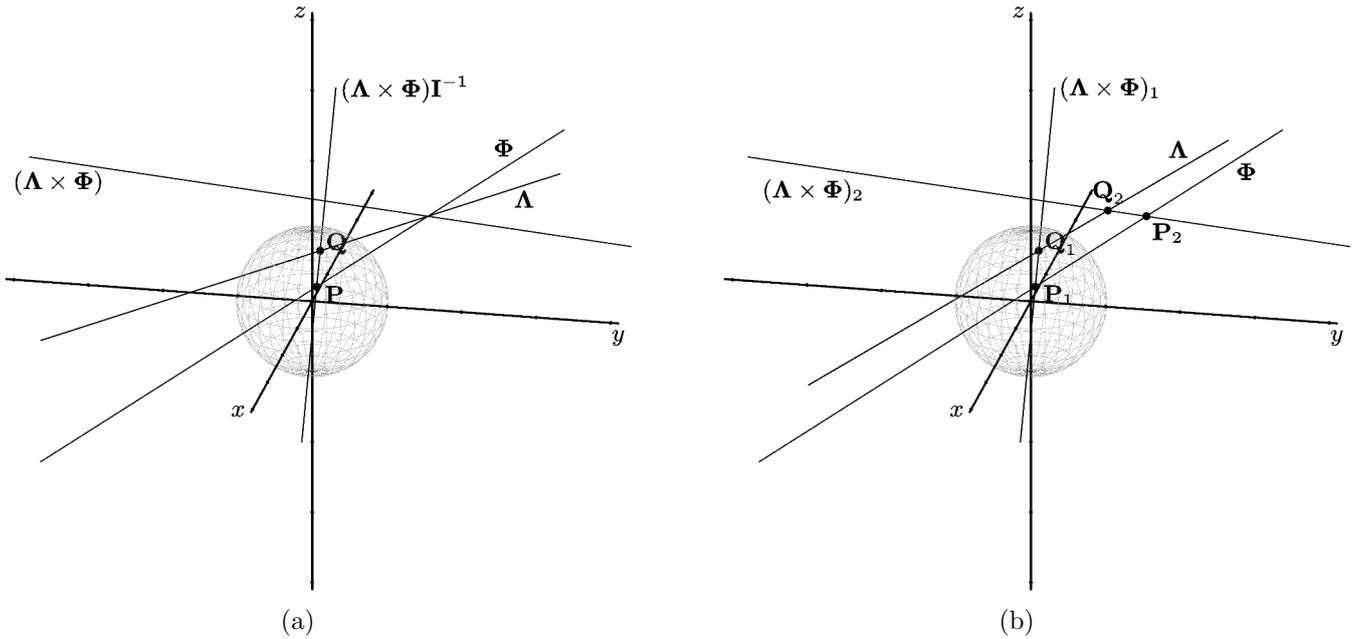
\hspace{-1cm}
\begin{subfloatenv}{ }
\begin{asy}
import Figure3D;
Figure f = Figure();
metric=Metric(Hyperbolic);

import solids;
triple f(pair t) { return (sin(t.x)*cos(t.y),sin(t.x)*sin(t.y),cos(t.x));}
surface sphere=surface(f,(0,0),(pi,2pi),16,32,Spline);
draw(surface(sphere),nullpen,meshpen=rgb(0.6,0.6,0.6),render(compression=Low,merge=true));

MV L = Line(-3/2,1,-1/2,-1,-5/2,-2); //Line(-3/2,1,-1,-1,-5/2,-1);
MV K = Line(1,5/3,-2,-1,3,2); //Line(1,4/3,-2,0,3,2);

MV F = cross(L,K);
MV F1=axis1(F);
MV F2=axis2(F);
MV Q1 = wedge(F1,join(L,origin)); Q1/=norm(Q1);
MV P1 = wedge(F1,join(K,origin)); P1/=norm(P1);
MV Q2 = wedge(F2,join(L,origin)); Q2/=norm(Q2);
MV P2 = wedge(F2,join(K,origin)); P2/=norm(P2);

MV M1 = join( Q1, P2 );
MV M2 = join( P1, P2 );

f.line(M1, "$\boldsymbol{\Lambda}$", position=0.9,align=(0,1,-0.5),draw_orientation=false);
f.line(M2, "$\boldsymbol{\Phi}$", position=0.1,align=(0,-0.5,1),draw_orientation=false);

MV G = cross(M1,M2);
f.line(G,"$(\boldsymbol{\Lambda}\times\boldsymbol{\Phi})$",draw_orientation=false,size=4,position=1,align=(0,1,-2));
f.line(G*I,"$(\boldsymbol{\Lambda}\times\boldsymbol{\Phi}){\boldsymbol{\textbf{I}}}^{-1}$",draw_orientation=false,size=4,position=0,align=(0,1,0));

MV Q = wedge(G*I,join(M1,origin)); Q/=norm(Q);
MV P = wedge(G*I,join(M2,origin)); P/=norm(P);

f.point(Q,"$\textbf{Q}$",align=(0,1,0.5),draw_orientation=false,draw_helper_lines=false);
f.point(P,"$\textbf{P}$",align=(0,1.25,-0.5),draw_orientation=false,draw_helper_lines=false);
asinh(norm(join(P,Q)));

\end{asy}
\end{subfloatenv}\hfill%
\begin{subfloatenv}{ }
\begin{asy}
import Figure3D;
Figure f = Figure();
metric=Metric(Hyperbolic);

import solids;
triple f(pair t) { return (sin(t.x)*cos(t.y),sin(t.x)*sin(t.y),cos(t.x));}
surface sphere=surface(f,(0,0),(pi,2pi),16,32,Spline);
draw(surface(sphere),nullpen,meshpen=rgb(0.6,0.6,0.6),render(compression=Low,merge=true));

MV L = Line(-3/2,1,-1/2,-1,-5/2,-2); //Line(-3/2,1,-1,-1,-5/2,-1);
MV K = Line(1,5/3,-2,-1,3,2); //Line(1,4/3,-2,0,3,2);

f.line(L, "$\boldsymbol{\Lambda}$", position=0.1,align=(0,-0.5,1),draw_orientation=false);
f.line(K, "$\boldsymbol{\Phi}$", position=0.9,align=(0,1,-0.5),draw_orientation=false);

MV F = cross(L,K);
MV F1=axis1(F);
MV F2=axis2(F);

f.line(F1,"$(\boldsymbol{\Lambda}\times\boldsymbol{\Phi})_1$",draw_orientation=false,size=4,position=1);
f.line(F2,"$(\boldsymbol{\Lambda}\times\boldsymbol{\Phi})_2$",position=0.8,align=(0,-1,-1),draw_orientation=false,size=4);

MV Q1 = wedge(F1,join(L,origin)); Q1/=norm(Q1);
MV P1 = wedge(F1,join(K,origin)); P1/=norm(P1);

f.point(Q1,"$\textbf{Q}_1$",align=(0,1,0.5),draw_orientation=false,draw_helper_lines=false);
f.point(P1,"$\textbf{P}_1$",align=(0,1.25,-0.5),draw_orientation=false,draw_helper_lines=false);
asinh(norm(join(P1,Q1)));

MV Q2 = wedge(F2,join(L,origin)); Q2/=norm(Q2);
MV P2 = wedge(F2,join(K,origin)); P2/=norm(P2);

f.point(Q2,"$\textbf{Q}_2$",align=(0,0,1),draw_orientation=false,draw_helper_lines=false);
f.point(P2,"$\textbf{P}_2$",align=(0,1,-1),draw_orientation=false,draw_helper_lines=false);
asinh(norm(join(P2,Q2)));

\end{asy}
\end{subfloatenv}
\caption{Commutator in \Hy{3}}
\label{commutator in H3}
\end{figure}

The commutator \(\mb{\Lambda}\times\mb{\Phi}\) of two intersecting proper lines \(\mb{\Lambda}\) and \(\mb{\Phi}\) is a simple bivector, i.e.\ it is a line.
The commutator is a proper, null, or improper line if the intersection point is proper, null, or improper, respectively.
If the proper lines intersect at a proper point, the commutator is perpendicular to both lines.
Proper lines intersecting at a null point are called parallel.
If proper lines intersect at an improper point, they are called hyperparallel.
The minimal separation between such lines is attained along the line \((\mb{\Lambda}\times\mb{\Phi})\I^{-1}\), 
whose polar line is the commutator product (see figure~\ref{commutator in H3}(a)).

If proper lines do not intersect, the commutator is a non-simple bivector, whose axes pass through the lines.
For example, see figure~\ref{commutator in H3}(b), where \((\mb{\Lambda}\times\mb{\Phi})_1\) is the proper axis of the commutator and \((\mb{\Lambda}\times\mb{\Phi})_2\) is its improper axis.
If the lines are proper and normalised, then the distance \(r\) and the angle \(\alpha\) between them obey 
\(\mb{\Lambda}\cdot\mb{\Phi}=-\cosh r\cos\alpha\) and \(|\mb{\Lambda}\vee\mb{\Phi}|=\sinh r\sin\alpha\),
which gives
\begin{equation}
2\sinh^2 r=u^2+v^2-1+\sqrt{(u^2+v^2-1)^2+4v^2},
\end{equation}
and \(\cos\alpha=-u/\cosh r\), where \(u=\mb{\Lambda}\cdot\mb{\Phi}\) and \(v=\mb{\Lambda}\vee\mb{\Phi}\).

The rotation of, say, \(\mb{\Lambda}\) by the angle \(\alpha\) around the proper axis of the commutator
brings \(\mb{\Lambda}\) in contact with \(\mb{\Phi}\) at an improper point without changing the distance between them.
As \(\mb{\Lambda}\) approaches \(\mb{\Phi}\) the axes of the commutator retain their attitude but not their weight.
At the point of contact, the commutator becomes simple: the proper axis disappears (compare figures~\ref{commutator in H3}(a) and~\ref{commutator in H3}(b)).
Similarly, a translation of \(\mb{\Lambda}\) by the distance \(r\) generated by the improper axis brings \(\mb{\Lambda}\)
in contact with \(\mb{\Phi}\) at a proper point without changing the angle between them.
At the point of contact, the commutator becomes simple: the improper axis disappears.

The usual formulas for projections, rejections, and reflections apply in \Hy{3}.
As long as the geometric objects involved in a projection and reflection are proper, the result of these two transformations is also a proper geometric object.
On the other hand, the result of a rejection of a proper geometric object in another proper geometric object is improper.
Moreover, the rejection falls on the polar counterpart of the object with respect to which the rejection is performed.
Projections, rejections, and reflections in null geometric objects are not possible.

The spinors are defined as usual and every spinor \(S\) can be written as \(S=e^A\) or \(S=-e^A\) where \(A\) is a bivector.
Any proper motion of a geometric object \(A_k\) can be obtained as an action \(SA_kS^{-1}\) of spinor \(S\).
Four different types of proper motion are possible in \Hy{3}: rotation around a proper line, translation along a proper line,
null translation, and a combination of rotation around a proper line and translation along the same line.

\begin{figure}[t]
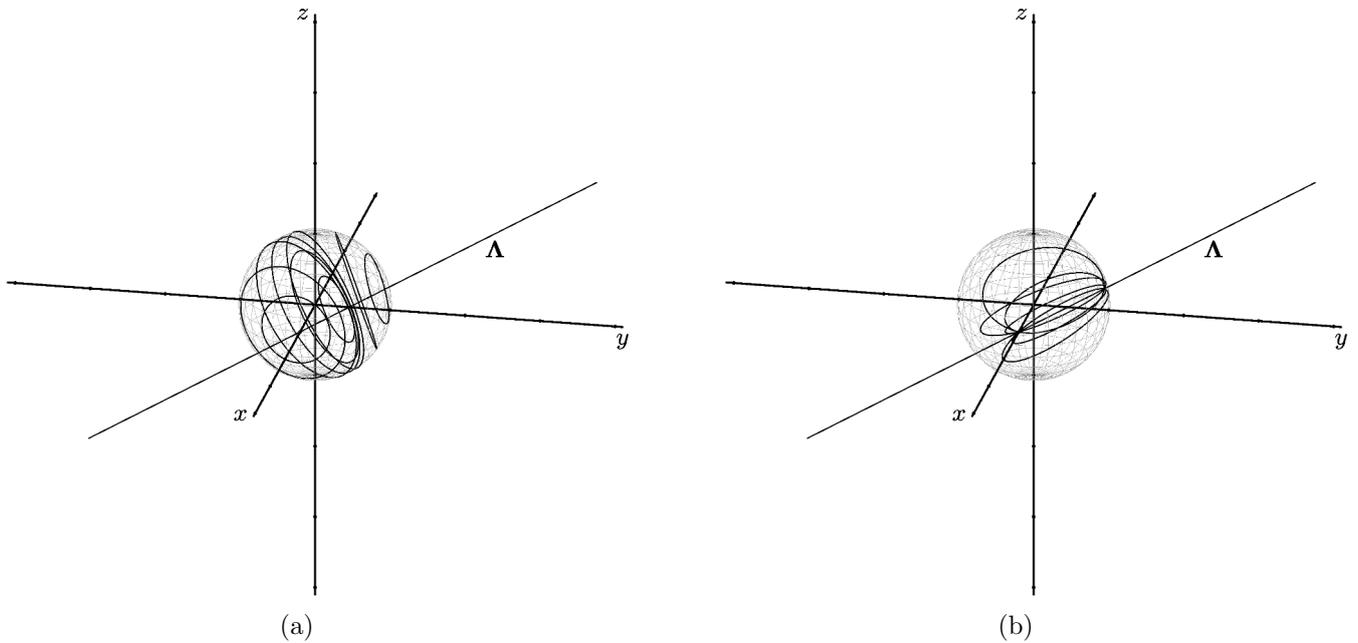
\hspace{-1cm}
\begin{subfloatenv}{ }
\begin{asy}
import Figure3D;
Figure f = Figure();
metric=Metric(Hyperbolic);

import solids;
triple f(pair t) { return (sin(t.x)*cos(t.y),sin(t.x)*sin(t.y),cos(t.x));}
surface sphere=surface(f,(0,0),(pi,2pi),16,32,Spline);
draw(surface(sphere),nullpen,meshpen=rgb(0.6,0.6,0.6),render(compression=Low,merge=true));

var P = Point(1,1,0,0);
var Q = Point(1,0,1,1/3);
var L = join(P,Q);
L/=norm(L);
write(L);

f.line(L, "$\boldsymbol{\Lambda}$", align=(0,0,-1),draw_orientation=false);

MV C = toline(L).centre();
MV K = join( C, origin);
real[] coef = { -2, -1, -1/2 };
MV[] points_centre = {};
for( real c: coef ) { MV S = exp( c*I*K ); MV M = S*C/S; points_centre.push(M); }

for(MV P: points_centre) { \
  triple f(real t) { MV S = cos(1/2*t)-L*sin(1/2*t); return totriple(S*P/S); }; \
  real x(real t) { return f(t).x; }; \
  real y(real t) { return f(t).y; }; \
  real z(real t) { return f(t).z; }; \
  path3 c = graph(x, y, z, -pi, pi, operator ..); \
  draw(c, currentpen+0.5); \
}

MV normalize( MV m ) { return m/norm(m);};
MV[] points_left = {};
for( MV P: points_centre ) { MV S = exp(I*normalize(dot(L,P)/P)); MV M = S*P/S; points_left.push(M); }

for(MV P: points_left) { \
  triple f(real t) { MV S = cos(1/2*t)-L*sin(1/2*t); return totriple(S*P/S); }; \
  real x(real t) { return f(t).x; }; \
  real y(real t) { return f(t).y; }; \
  real z(real t) { return f(t).z; }; \
  path3 c = graph(x, y, z, -pi, pi, operator ..); \
  draw(c, currentpen+0.5); \
}

MV[] points_right = {};
for( MV P: points_centre ) { MV S = exp(-I*normalize(dot(L,P)/P)); MV M = S*P/S; points_right.push(M); }

for(MV P: points_right) { \
  triple f(real t) { MV S = cos(1/2*t)-L*sin(1/2*t); return totriple(S*P/S); }; \
  real x(real t) { return f(t).x; }; \
  real y(real t) { return f(t).y; }; \
  real z(real t) { return f(t).z; }; \
  path3 c = graph(x, y, z, -pi, pi, operator ..); \
  draw(c, currentpen+0.5); \
}

\end{asy}
\end{subfloatenv}\hfill%
\begin{subfloatenv}{ }
\begin{asy}
import Figure3D;
Figure f = Figure();
metric=Metric(Hyperbolic);

import solids;
triple f(pair t) { return (sin(t.x)*cos(t.y),sin(t.x)*sin(t.y),cos(t.x));}
surface sphere=surface(f,(0,0),(pi,2pi),16,32,Spline);
draw(surface(sphere),nullpen,meshpen=rgb(0.6,0.6,0.6),render(compression=Low,merge=true));

var P = Point(1,1,0,0);
var Q = Point(1,0,1,1/3);
var L = join(P,Q);
L/=norm(L);
write(L);

f.line(L, "$\boldsymbol{\Lambda}$", align=(0,0,-1),draw_orientation=false);

MV C = toline(L).centre();
MV K = join( C, origin);
real[] coef = { -1, -1/2, 1/2 };
MV[] points_centre = {};
for( real c: coef ) { MV S = exp( c*I*K ); MV M = S*C/S; points_centre.push(M); }

for(MV P: points_centre) { \
  triple f(real t) { MV S = cosh(1/2*t) - I*L*sinh(1/2*t); return totriple(S*P/S); }; \
  real x(real t) { return f(t).x; }; \
  real y(real t) { return f(t).y; }; \
  real z(real t) { return f(t).z; }; \
  path3 c = graph(x, y, z, -10, 10, operator ..); \
  draw(c, currentpen+0.5); \
}

MV[] points_up = {};
for( MV P: points_centre ) { MV S = exp(-1/2*0.5*L); MV M = S*P/S; points_up.push(M); }

for(MV P: points_up) { \
  triple f(real t) { MV S = cosh(1/2*t) - I*L*sinh(1/2*t); return totriple(S*P/S); }; \
  real x(real t) { return f(t).x; }; \
  real y(real t) { return f(t).y; }; \
  real z(real t) { return f(t).z; }; \
  path3 c = graph(x, y, z, -10, 10, operator ..); \
  draw(c, currentpen+0.5); \
}

MV[] points_down = {};
for( MV P: points_centre ) { MV S = exp(1/2*0.5*L); MV M = S*P/S; points_down.push(M); }

for(MV P: points_down) { \
  triple f(real t) { MV S = cosh(1/2*t) - I*L*sinh(1/2*t); return totriple(S*P/S); }; \
  real x(real t) { return f(t).x; }; \
  real y(real t) { return f(t).y; }; \
  real z(real t) { return f(t).z; }; \
  path3 c = graph(x, y, z, -10, 10, operator ..); \
  draw(c, currentpen+0.5); \
}

\end{asy}
\end{subfloatenv}
\caption{Rotation and translation in \Hy{3}}
\label{rotation and translation in H3}
\end{figure}

The action of a spinor
\begin{equation}
S=e^{-\tfrac{1}{2}(\alpha+\lambda\I)\mb{\Lambda}},
\end{equation}
where \(\mb{\Lambda}\) is a normalized proper line, yields a combination of the rotation around \(\mb{\Lambda}\) by the angle \(\alpha\)
and the translation along the same line by distance \(|\lambda|\) (the direction of the translation depends on the sign of \(\lambda\)).
Pure rotation obtains if \(\lambda=0\) and pure translation obtains if \(\alpha=0\).
An example of the rotation around line \(\mb{\Lambda}\) is given in Figure~\ref{rotation and translation in H3}(a), 
where the circular trajectories (the rotation angle varies from \(-\pi\) to \(\pi\)) of nine different points are shown.
The translation along the same line (or equivalently the translation generated by the line polar to \(\mb{\Lambda}\))
is given in Figure~\ref{rotation and translation in H3}(b) for nine other points (the translation distance varies from \(-\infty\) to \(+\infty\)).
All points on a given trajectory are at a constant distance from  \(\mb{\Lambda}\).

The action of a spinor
\begin{equation}
S=e^{-\tfrac{1}{2}\theta\mb{\Lambda}},
\end{equation}
where \(\mb{\Lambda}\) is a null line, yields a null translation anchored at the point where \(\mb{\Lambda}\) touches the unit sphere.
The interpretation of \(\theta\) depends on the weight of \(\mb{\Lambda}\).

%\begin{figure}[t]\hspace{-1cm}
%\begin{subfloatenv}{ }
%\begin{asy}
%import Figure3D;
%Figure f = Figure();
%metric=Metric(Hyperbolic);
%
%import solids;
%triple f(pair t) { return (sin(t.x)*cos(t.y),sin(t.x)*sin(t.y),cos(t.x));}
%surface sphere=surface(f,(0,0),(pi,2pi),16,32,Spline);
%draw(surface(sphere),nullpen,meshpen=rgb(0.6,0.6,0.6),render(compression=Low,merge=true));
%
%var P = Point(1,1,0,0);
%var Q = Point(1,0,1,1/3);
%var L = join(P,Q);
%MV M = 7*L+I*L;
%M/=norm(M);
%write(M);
%
%f.line(L, "$\boldsymbol{\Lambda}$", align=(0,0,-1),draw_orientation=false);
%
%triple f(real t) { MV S = exp(-1/2*t*M); MV Q = S*origin/S; if ( Q.is_point_excess() > 1e-2 ) { return (0,0,0); }; return totriple(Q); }; \
%real x(real t) { return f(t).x; }; \
%real y(real t) { return f(t).y; }; \
%real z(real t) { return f(t).z; }; \
%path3 c = graph(x, y, z, -40, 0, operator ..); \
%draw(c, currentpen+0.5); \
%
%\end{asy}
%\end{subfloatenv}\hfill%
%\begin{subfloatenv}{ }
%\end{subfloatenv}
%\caption{Combined rotation and translation in \Hy{3}}
%\label{combined rotation and translation in H3}
%\end{figure}

\bibliographystyle{plain}
\bibliography{g}

\end{document}